\def\Draft{0}

\ifnum\Draft=1
	\documentclass[12pt,onecolumn,draftclsnofoot]{IEEEtran}
\else
	\documentclass[]{IEEEtran}
\fi

\usepackage{graphicx,xcolor,tikz} % graphic packages
\usepackage{amssymb,amsmath,mathptmx,mathrsfs} % math packages
\usepackage[algoruled,vlined]{algorithm2e} % algorithm package
\usepackage{hyperref} % formating packages
\usepackage[nocompress]{cite} % citations in separated brackets and sorted by number

\SetKwInput{KwPara}{Parameters} % Input type "Parameters"

\newcommand{\R}{\mathbb{R}}
\newcommand{\C}{\mathbb{C}}
\newcommand{\N}{\mathbb{N}}
\newcommand{\F}{\mathcal{F}}
\newcommand{\G}{\mathcal{G}}

\DeclareMathOperator*{\argmin}{arg\,min}
\DeclareMathOperator*{\argmax}{arg\,max}
\DeclareMathOperator*{\supp}{supp}

\renewcommand{\u}[1]{\underline{#1}\hspace*{1pt}}
\newcommand{\n}[1]{\widetilde{#1}}

\newcommand{\domM}{\mathscr{M}} % spatial domain
\newcommand{\domP}{\mathscr{P}} % parameter domain

\newcommand{\ptM}{m} % point in spatial domain
\newcommand{\ptP}{p} % point in parameter domain

\newcommand{\setI}{I} % index set (support)
\newcommand{\iS}{J} % index set (support choice)
\newcommand{\setC}{K} % index set (neighborhood)

\newcommand{\nbrM}{N_2} % number of sampling points in spatial domain
\newcommand{\nbrT}{N_1} % number of sampling points in time domain
\newcommand{\nbrP}{N_3} % number of sampling points in parameter domain

\newcommand{\dicM}{\u{D}} % dictionary matrix
\newcommand{\splM}{\u{B}} % sample matrix
\newcommand{\splV}{\u{b}} % sample vector

\newcommand{\resM}{\u{R}} % residuum matrix
\newcommand{\resV}{\u{r}} % residuum vector

\newcommand{\solM}{\u{X}} % solution matrix
\newcommand{\solV}{\u{x}} % solution vector

\newcommand{\fS}{\mathbb{J}} % feasible support set
\newcommand{\mtx}{\u{M}} % suppport matrix of feasible set
\newcommand{\p}{\alpha} % Parameter for support constraint
\newcommand{\pp}{\gamma} % Parameter for Lipschitz constraint

\newtheorem{thm}{Theorem}
\newtheorem{defi}{Definition}
\newtheorem{rmk}{Remark}

\begin{document}

\title{Generalized orthogonal matching pursuit for multiple measurements - A structural approach}

\author{Florian Bo{\ss}mann\IEEEauthorrefmark{1}\thanks{\IEEEauthorrefmark{1}University of Passau, Mathmatics / Digital Image Processing, Florian.Bossmann@uni-passau.de}}

\maketitle

\begin{abstract}
Sparse data approximation has become a popular research topic in signal processing. However, in most cases only a single measurement vector (SMV) is considered. In applications, the multiple measurement vector (MMV) case is more usual, i.e., the sparse approximation problem has to be solved for several data vectors coming from closely related measurements. Thus, there is an unknown inter-vector correlation between the data vectors. Using SMV methods typically does not return the best approximation result as the correlation is ignored. In the past few years several algorithms for the MMV case have been designed to overcome this problem. Most of these techniques focus on the approximation quality while quite strong assumptions to the type of inter-vector correlation are made.

While we still want to find a sparse approximation, our focus lies on preserving (possibly complex) structures in the data. Structural knowledge is of interest in many applications. It can give information about e.g., type, form, number or size of objects of interest. This may even be more useful than information given by the non-zero amplitudes itself. Moreover, it allows efficient post processing of the data. We numerically compare our new approach with other techniques and demonstrate its benefits in two applications.
\end{abstract}


\section{Introduction}

Sparse approximations of given data are of great interest in many different applications. They are used in image processing for e.g., denoising \cite{Elad06,Dadov07}, compression \cite{LePennec05} or restoration \cite{Mairal08}. Sparsity assumptions appear in face and speech recognition \cite{Katz87,Wright09}, magnetic resonance imaging (MRI) and computer tomography (CT) \cite{Lustig07,Sidky08}, as well as in non-destructive testing \cite{Kaaresen99,Bossmann12} and seismic data processing \cite{Scales88,Trad03,Bossmann15,Bossmann16}. A detailed overview can also be found in \cite{Starck10,Rao98} and the references therein.

The sparse approximation itself can be stated as follows: Given a dictionary matrix $\dicM\in\C^{\nbrT\times\nbrP}$ and a measurement vector $\splV\in\C^{\nbrT}$, solve
\begin{align}\label{eq:repProb}
\min\limits_{\solV\in\C^{\nbrP}}\|\solV\|_0 \text{ s.t. }\|\dicM\solV-\splV\|_2\leq\varepsilon
\end{align}
for a given $\varepsilon>0$. Here $\|\solV\|_0$ is the $\ell_0$-quasi-norm, i.e., $\|\solV\|_0:=|\{k\ |\ x_k\neq0\}|$. The vector $\solV$ is said to be $L$-sparse, if $\|\solV\|_0\leq L$ for $L\in\N$. The matrix $\dicM$ is constructed using a basis, frame or dictionary in which the given data $\splV$ is assumed to be sparse. Typical examples are Fourier \cite{Candes06} or Wavelet bases \cite{Mallat99} as well as Curvelet \cite{Candes00} or Shearlet frames \cite{Labate05}. Dictionaries can be designed according to the underlying application, as e.g., the Gabor impulse in ultrasonic testing \cite{Demirli01,Bossmann12}.

The exact solution of (\ref{eq:repProb}) can in general only be found combinatorially, i.e., by considering all possible supports of $\solV$. Hence, finding the exact solution becomes NP-hard. There are two main strategies to find at least an approximate solution of (\ref{eq:repProb}): The first strategy is, replacing the $\ell_0$-quasi-norm by the $\ell_1$-norm what makes the problem convex. This approach is known as convex relaxation or basis pursuit \cite{Chen01,Tropp04}. Greedy algorithms are another strategy to solve (\ref{eq:repProb}) approximatively. Those methods iteratively built up a global approximation by solving local subproblems \cite{Temlyakov11}. Matching Pursuit (MP) and Orthogonal Matching Pursuit (OMP) \cite{Tropp04} may be the most known algorithms in this context. In recent years, more advanced algorithms have been developed such as Stagewise OMP \cite{Donoho12}, Compressive Sampling MP \cite{Needell08,Needell09} and regularized OMP \cite{Needell08,Needell09a}. An overview can be found in \cite{Tropp10}.

Eq. (\ref{eq:repProb}) is known as single measurement vector problem (SMV). It has been studied extensively over the last few years. However, there is an extension known as multiple measurement vector problem (MMV). Instead of only having one data vector $\splV\in\C^{\nbrT}$, several measurements $\splM:=(\splV^1,\ldots,\splV^{\nbrM})\in\C^{\nbrT\times\nbrM}$ are given. The problem is stated similar as
\begin{align}\label{eq:repProbMMV}
\min\limits_{\solM\in\C^{\nbrP\times\nbrM}}\|\solM\|_{0,\infty} \text{ s.t. }\|\dicM\solM-\splM\|_2\leq\varepsilon
\end{align}
with $\solM=(\solV^1,\ldots,\solV^{\nbrM})$ and $\|\solM\|_{0,\infty}=\max_k\|\solV^k\|_0$, i.e., each vector $\solV^k$ is sparse. In fact this problem seems to be more common in many applications. It appears e.g., in non-destructive testing \cite{Kaaresen99,Bossmann12}, in seismic data \cite{Scales88,Trad03,Bossmann15,Bossmann16} or Magnetoencephalography (MEG) \cite{Cotter05}. The MMV formulation can be used, whenever several measurements of the same or similar objects were made. One straight-forward approach to solve (\ref{eq:repProbMMV}) is, to split it into $\nbrM$ SMV problems and solve these independently using the methods mentioned above. However, our intuition tells us, since the $\nbrM$ measurements were made in a quite similar set-up, also the obtained data vectors $\splV^1,\ldots,\splV^{\nbrM}$ should be correlated somehow. Simply solving $\nbrM$ SMV problems ignores this correlation and thus the solution quality might suffer.

Recently, new methods have been developed that consider inter-signal correlation. In \cite{Cotter05} an extension for MP and the FOCal Underdetermined System Solver (FOCUSS) are presented. Bayesian methods are considered in \cite{Zhang10,Ziniel13}. In \cite{Tropp06a,Tropp06b} the authors introduce Greedy pursuit and convex relaxation for the MMV problem. Theoretical results have been shown e.g., in \cite{Chen06}. All these methods force a common support in the reconstructed solution, i.e., the reconstructed matrix $\solM$ has only few nonzero rows or, in other words, the columns of $\solM$ have (nearly) the same support. In \cite{Baron06} two joint sparsity models (JSM) for compressed sensing are introduced. JSM-1 considers solutions where all columns $\solV^k$ can be written as the sum $\solV^k=\solV^{c}+\solV^{k,u}$ of a common sparse component $\solV^c$ that is equal for each column and another unique sparse vector $\solV^{k,u}$. JSM-2 is equal to the support constraint considered above. Another approach is presented in \cite{Zheng11,Lu12} where correlated measurements are assumed to have sparse approximations that are close in the euclidean distance. This idea is related to dynamic compressed sensing \cite{Ziniel13a,Angelosante09}. Here neighboring columns are assumed to have similar support. In both cases, the support is allowed to change slowly over different data vectors. In most cases the used methods penalize non-smooth rows in $\solM$.

However, all methods have quite restricting support assumptions and hence cannot reconstruct simple geometries in the solution. As an example consider $\solM$ to be the identity matrix. There is no common support between all columns and the rows are non-smooth. Nevertheless, the matrix is still clearly structured. The linear structure can be described by a shift of $1$ index per column. In the next section we introduce a generalized version of orthogonal matching pursuit for multiple measurements (GM-OMP) that takes complex structures in the data into account. Numerical evidence for the proposed method are shown in the third section of this work.

\section{The Algorithm}

In this section we first introduce OMP and discuss its generalization to the MMV problem. GM-OMP increases the support of the solution $\solM$ in each iteration by adding an index set $\iS\in\fS$ where $\fS$ is the set of feasible selections. The parametrization and selection of $\iS$ is the main idea of GM-OMP and is discussed in the second subsection. After first theoretical results are shown in the third subsection, we present an a-posteriori denoising technique that is based on the structural component of the reconstructed solution.

\subsection{OMP and GM-OMP}

Orthogonal matching pursuit is a greedy algorithm that seeks to find a sparse solution of (\ref{eq:repProb}). For simplicity we assume that the columns of $\dicM$ are normalized. Then the iterative scheme of OMP can be summarized as follows:
\begin{enumerate}
\item Set the residual $\resV=\splV$ and the support $\setI=\emptyset$.
\item Calculate $i=\argmax|\dicM^*\resV|$ and update $\setI\leftarrow\setI\cup\{i\}$.
\item Solve $\solV=\argmin\limits_{\supp \underline y\subseteq\setI}\|\splV-\dicM \underline y\|_2$ and set $\resV=\splV-\dicM\solV$.
\item Iterate 2-3 until a stopping criterion holds.
\end{enumerate}
Here $\dicM^*$ is the transposed conjugate complex matrix. Hence, the algorithm chooses the column of $\dicM$ that correlates most with the residual and adds its index to the support set in step 2. Step 3 calculates the best approximation according to the selected support. The algorithm may e.g., be stopped after $L$ iterations (the solution $\solV$ is $L$-sparse then), or when the residuum drops below a threshold, i.e., $\|\resV\|_2\leq\varepsilon$.

Now, let us consider the MMV problem shown in (\ref{eq:repProbMMV}). The idea of GM-OMP is surprisingly simple. We only adapt the second step of OMP, while all other steps stay the same. Therefore, note that OMP chooses one index $i$ and adds it to the support set $\setI$. Since we are now dealing with multiple measurements, GM-OMP is allowed to add not only one index $i$, but multiple indices to the support (e.g., one index per column of $\solM$). Let us denote the set of all indices added by $\iS$. This index set should be chosen from a set of feasible selections $\iS\in\fS\subseteq\mathcal{P}(\{(j,k)\ | j\leq\nbrP,k\leq\nbrM\})$ where $\mathcal{P}(\cdot)$ denotes the power set. The second step of GM-OMP now reads as follows:
\begin{enumerate}
\setcounter{enumi}{1}
\item Choose $\iS\in\fS$ and update $\setI\leftarrow\setI\cup\iS$.
\end{enumerate}
The complete scheme of GM-OMP is shown in Alg. \ref{alg:GMOMP} where the maximum number of iterations $L$ and the minimal residuum norm $\varepsilon_R$ is included as stopping criterion. Of course we need to define the feasible set $\fS$ and find a suitable choice $\iS\in\fS$. This problem will be discussed in the next subsection. Let us first consider three examples to clarify the principle of GM-OMP and the feasible set $\fS$.

\begin{algorithm}
  \KwData{$\splM,\dicM$}
  \KwPara{$\varepsilon_R,L$}
  Set $\resM=\splM$ and $\setI=\emptyset$\;
  \For{$l=1,\ldots,L$}{
  	Choose $\iS\in\fS$ and update $\setI\leftarrow\setI\cup\iS$\;
  	Solve $\solM=\argmin\limits_{\supp \underline Y\subseteq\setI}\|\splM-\dicM \underline Y\|_2$ and set $\resM=\splM-\dicM\solM$\;
  	\lIf{$\|\resM\|_2\leq\varepsilon_R$}{stop}
  }
  \caption{GM-OMP}
  \label{alg:GMOMP}
\end{algorithm}

For our first example, consider the set
\begin{align}\label{def:POMPset}
\fS=\fS_P=\{\supp\mtx\ | \text{ each column of }\mtx\text{ is 1-sparse}\},
\end{align}
i.e., $\fS_P$ contains all index sets that can be associated with the support of matrices $\mtx$ having at most one non-zero element per column. Here $\mtx$ is of same size as $\solM$. Choose $\iS\in\fS_P$ such that
\begin{align*}
(j,k)\in\iS\ \Leftrightarrow j=\argmax|\dicM^*\resV^k|,
\end{align*}
where $\resV^k$ is the $k$-th column of the residual matrix $\resM$. This way, GM-OMP is equivalent to OMP parallelly applied to each column of $\splM$. In our next example define
\begin{align}\label{def:VOMPset}
\fS=\fS_V=\{\supp\mtx\ |\ \|\mtx\|_{0}=1\}
\end{align}
as the set of all 1-sparse supports. It is easy to see that the best choice $\iS\in\fS_V$ is $\iS=\{\argmax|\dicM^*\resM|\}$ the index set containing only the position of the maximum absolute value of $\dicM^*\resM$. Now, GM-OMP is identical to using OMP on the vectorized formulation of (\ref{eq:repProbMMV}), i.e., rewrite $\solM,\splM$ as column vectors and $\dicM$ becomes a block diagonal matrix. As our last example, consider
\begin{align}\label{def:SOMPset}
\fS=\fS_S=\{\supp\mtx\ |\ \mtx\text{ has at most one non-zero row }\}
\end{align}
containing all support sets with constant row index. A possible choice $\iS\in\fS_S$ is given by
\begin{align}\label{eq:choiceSOMP}
i=\argmax_j\|(\dicM^*\resM)_{j,\cdot}\|_{\lambda}
\end{align}
and $\iS=\{(i,j)\ |\ j\leq\nbrM\}$. Here $\|(\dicM^*\resM)_{j,\cdot}\|_{\lambda}$ denotes the $\lambda$-norm of the $j$-th row of the matrix. For $\lambda=1$ GM-OMP becomes the simultaneous OMP (S-OMP) introduced in \cite{Tropp06a}, the case $\lambda=2$ is discussed in \cite{Cotter05}.

For the three demonstrated choices of $\fS$, GM-OMP transforms into well known algorithms for the MMV problem. However, neither $\fS_P$ nor $\fS_V$ contain structured sets while $\fS_S$ is bounded to row-sparsity of $\solM$ (see discussion in the introduction). Thus, we introduce a more general choice for $\fS$ in the next subsection. Here, $\fS=\fS(\p,\pp)$ can be adapted by parameters. We will see that the three examples form the extreme cases of the parameter choice.

\subsection{Feasible set and selection}

Alg. \ref{alg:GMOMP} demonstrates the generalized scheme of OMP for multiple measurements. The set $\fS$ represents the feasible sparsity patterns that can be chosen per iteration. However, since $\fS$ is a subset of a power set, it can be of exponential size. Thus, it is not sufficient to leave it to the user as input data. Instead we will parametrize $\fS$ and define a selection rule for $\iS\in\fS$ based on the parameters. This way, the user only has to choose parameters that describe sparsity patterns suitable for his application.

\begin{rmk}
Following, we describe a parametrization idea that, in the authors opinion, can be used in many applications. However, the reader might choose a different description of $\fS$ and an according selection rule $\iS\in\fS$ that is more suitable for the particular problem.
\end{rmk}

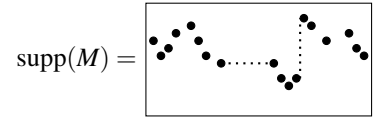
\begin{figure}
\centering
\begin{tikzpicture}
	\node[left] at (0,0.75) {$\supp(M)=$};
	\draw (0,0) rectangle (3,1.5);
	\draw[fill=black] (0.1,1) circle (0.05);
	\draw[fill=black] (0.2,0.8) circle (0.05);
	\draw[fill=black] (0.3,0.9) circle (0.05);
	\draw[fill=black] (0.4,1.1) circle (0.05);
	\draw[fill=black] (0.6,1.2) circle (0.05);
	\draw[fill=black] (0.7,1) circle (0.05);
	\draw[fill=black] (0.8,0.8) circle (0.05);
	\draw[fill=black] (1,0.7) circle (0.05);
	\draw[dotted,thick] (1,0.7) -- (1.7,0.7);
	\draw[fill=black] (1.7,0.7) circle (0.05);
	\draw[fill=black] (1.8,0.5) circle (0.05);
	\draw[fill=black] (1.9,0.4) circle (0.05);
	\draw[fill=black] (2,0.5) circle (0.05);
	\draw[dotted,thick] (2.05,0.5) -- (2.05,1.3);
	\draw[fill=black] (2.1,1.3) circle (0.05);
	\draw[fill=black] (2.2,1.2) circle (0.05);
	\draw[fill=black] (2.4,1) circle (0.05);
	\draw[fill=black] (2.7,1.1) circle (0.05);
	\draw[fill=black] (2.8,0.9) circle (0.05);
	\draw[fill=black] (2.9,0.8) circle (0.05);
\end{tikzpicture}
\caption{$\supp(M)$ as samples of a function possibly having a disconnected domain and discontinuities (dotted lines).}
\label{fig:suppM}
\end{figure}

Note that $\iS\in\fS$ is a set of two-dimensional elements (row and column indices). Our idea is, to use exactly two parameters $\p,\pp$ to determine $\fS=\fS(\p,\pp)$. It is clear that we cannot cover all sets $\fS$ with two parameters, since the number of possible choices for $\fS$ grows exponential. Thus, we need a parametrization that generates suitable sets for applications. Analogous to the given examples (\ref{def:POMPset})-(\ref{def:SOMPset}), we identify an element $\iS\in\fS$ by its pattern matrix $\mtx$ where $\iS=\supp\mtx$ holds. Since we assume the columns of $\solM$ to be sparse, it is reasonable to permit only matrices $\mtx$ with (at most) 1-sparse columns, i.e., in each iteration of GM-OMP the support of $\solM$ should at most grow by one index per column. Given such a matrix we interpret its sparsity pattern as samples of a function, mapping the column indices to corresponding row indices (Fig. \ref{fig:suppM}). Due to the 1-sparse columns of $\mtx$ this mapping is unique but not necessarily defined for all rows (there may be zero columns in $\mtx$). Having this in mind, we postulate the sparsity pattern of $\mtx$ to hold two conditions:
\begin{itemize}
\item The domain of the sparsity pattern should be connected.
\item The sparsity pattern should be (Lipschitz-)continuous.
\end{itemize}
Fig. \ref{fig:suppM} shows a sparsity pattern where both conditions do not hold (see the dotted lines). Next, we formulate our parametrization of $\fS$ that uses two parameters $\p,\pp$ to ensure the above stated conditions. Therefore, we introduce the parameter and measurement space.

Let $\domP$ and $\domM$ be metric spaces. For $\solM\in\C^{\nbrP\times\nbrM}$ let $\ptP_1,\ldots,\ptP_{\nbrP}\in\domP$ and $\ptM_1,\ldots,\ptM_{\nbrM}\in\domM$ be given. We call $\domP$ and $\domM$ the parameter space and measurement space respectively. The elements $\ptP_j$ are parameters (of the dictionary) and $\ptM_j$ is a measurement (setup).

\begin{rmk}
At first glance it seems quite restricting to require the existence of such spaces and elements. However, they come quite naturally. For example consider $\dicM$ being the Fourier matrix. Then $\ptP_j$ is the frequency of the $j$-th column of $\dicM$. If we use a Wavelet dictionary $\dicM$, the parameters $\ptP_j$ contain the shift and scaling of each column. For convolution matrices $\dicM$ each $\ptP_j$ is the shift of the $j$-th column. On the other hand, consider the measurement data $\splM$ was obtained using several sensors at different positions, each column of $\splM$ corresponding to one sensor. Thus we can set $\ptM_j$ to the position of the $j$-th sensor. While these parameters give the following formulas a more reasonable interpretation, one can surely just use $\ptP_j=j$ and $\ptM_j=j$.
\end{rmk}

Now we can formulate the above stated conditions on $\iS=\supp\mtx$. Using the points $\ptM_i$ as vertices and defining edges using the metric $d_{\domM}(\cdot,\cdot)$ defined on $\domM$, we can state the connected sparsity pattern condition as, the graph
\begin{align}\label{eq:condSupp}
\left(\{\ptM_i\}_{(i,j)\in\iS}\ |\ \{\overline{\ptM_i\ptM_{i'}},\ d_{\domM}(\ptM_i,\ptM_{i'})\leq\p\}\right)
\text{ is connected.}
\end{align}
%\max\limits_{(i,j)\in\iS}\min\limits_{(i',j')\in\iS,i'\neq i}d_{\domM}(\ptM_i,\ptM_{i'})\leq\p,
Lipschitz continuity of the pattern is ensured if
\begin{align}\label{eq:condLipp}
d_{\domP}(\ptP_j,\ptP_{j'})\leq\pp\ d_{\domM}(\ptM_i,\ptM_{i'}) && \forall(i,j),(i',j')\in\iS
\end{align}
with the metric $d_{\domP}$ on $\domP$. We define $\fS=\fS(\p,\pp)$ by
\begin{align}\label{eq:defFS}
\fS(\p,\pp)=\{\iS\ |\ \iS\text{ satisfies (\ref{eq:condSupp}) and (\ref{eq:condLipp})}\}.
\end{align}
By (\ref{eq:condLipp}) we also ensure 1-sparse columns of $\mtx$ ($\iS=\supp\mtx$) if $\pp<\infty$. For $\pp=\infty$ we use the convention $\infty\cdot0=\lim_{\pp\rightarrow\infty}\pp\cdot0=0$. We obtain the relations $\fS(\infty,\infty)=\fS_P$, $\fS(0,\pp)=\fS_V$ and $\fS(\infty,0)=\fS_S$, i.e., our parametrization covers the shown examples.

Given the set $\fS(\p,\pp)$ we need to choose $\iS\in\fS$. Like in OMP, we search for the support $\iS$ that maximizes the correlation between dictionary and residuum, i.e., we would like to solve
\begin{align}\label{eq:greedy}
\iS=\argmax\limits_{\iS'\in\fS(\p,\pp)}\|(\dicM^*\resM)_{\iS'}\|_{\lambda}
\end{align}
for some $\lambda\geq1$ (compare Eq. \ref{eq:choiceSOMP}). Intuitively, we want all values $(\dicM^*\resM)_{\iS}$ to be of same order of magnitude, i.e., $\lambda=1$ (or $\lambda=2$) may be a good choice. Unfortunately, we can state the following theorem which is proven in the next subsection:

\begin{thm}\label{thm:NPhard}
For $\lambda<\infty$ and arbitrary $\p,\pp$ problem (\ref{eq:greedy}) is NP-hard.
\end{thm}

Hence, we use a greedy algorithm to approximatively solve (\ref{eq:greedy}). Indeed, the algorithm returns the exact solution of (\ref{eq:greedy}) with $\lambda=\infty$. Starting with the correlation matrix $C=|\dicM^*\resM|$ we iteratively built a matrix $\mtx$ such that $\iS=\supp\mtx\in\fS$. Beginning with $\mtx=0$ we add the position of the maximum value of $C$ to the support, i.e., we calculate $(i,j)=\argmax (C)_{i',j'}$ and update $(\mtx)_{i,j}=1$. To ensure that (\ref{eq:condSupp}) is not penalized, we restrict ourself to indices $i'\in\setC$ where $\setC$ is the set of all indices for which (\ref{eq:condSupp}) holds. Afterwards, the chosen element $(C)_{i,j}$ and all elements in $C$ that violate (\ref{eq:condLipp}) are set to zero, This way, it is guaranteed that (\ref{eq:condLipp}) is fulfilled. The scheme is shown in Alg. \ref{alg:GMOMPchoice}. 

\begin{algorithm}
  \KwData{$C=|\dicM^*\resM|$}
  \KwPara{$\p,\pp$}
  Set $\mtx=0\in\{0,1\}^{\nbrP\times\nbrM}$, $\setC=\{1,\ldots,\nbrM\}$\;
  \While{$(C)_{K,\cdot}\neq0$}{
  	$(i,j)=\argmax (C)_{i',j'}$ s.t. $i'\in\setC$, $j\leq\nbrP$\;
  	$(\mtx)_{i,j}=1$, $(C)_{i,j}=0$\;
  	$(C)_{i',j'}=0$ where $\supp\mtx\cup\{(i',j')\}$ violates (\ref{eq:condLipp})\;
  	$\setC=\{i\ |\ \min\limits_{(i',j')\in\supp\mtx}d_{\domM}(\ptM_i,\ptM_{i'})\leq\p,\ i\leq\nbrM\}$\;
  }
  \caption{GM-OMP greedy choice}
  \label{alg:GMOMPchoice}
\end{algorithm}

For implementation we recommend to replace $(C)_{I,\cdot}\neq0$ by $(C)_{I,\cdot}>\varepsilon$ using a reasonable threshold $\varepsilon$ (see also the discussion in the theory part). Furthermore, if $C=|\dicM^*\resM|$ from start on has zero entries (or elements below the threshold), they will never be chosen by the algorithm. This assures that the support is not artificially enlarged, i.e., $\supp\mtx\subseteq\supp C$ always holds.

It is easy to see that Alg. \ref{alg:GMOMPchoice} returns the exact solution for $(\p,\pp)=(\infty,\infty)$ and $(\p,\pp)=(0,\pp)$ (i.e,, for $\fS_P$ and $\fS_V$). Setting $(\p,\pp)=(\infty,0)$ (i.e., $\fS_S$) Alg. \ref{alg:GMOMPchoice} solves (\ref{eq:choiceSOMP}) for $\lambda=\infty$.

\subsection{Theoretical results}

We first prove Theorem \ref{thm:NPhard} that is the motivation for Alg. \ref{alg:GMOMPchoice}.

\begin{IEEEproof}[Proof of Theorem \ref{thm:NPhard}]
We give a polynomial-time reduction of the coloring problem: Given a graph $G$ and a number of colors $C$, assign a color to each vertex $v_j\in G$ such that for each edge $\overline{v_iv_j}\in G$ the vertices $v_i,v_j$ have a different color.

Let $G$ have $M$ vertices, each vertex with at most $n$ edges. For the reduction we need each vertex to have the same amount of edges. Thus, we add edges $\overline{v_j\cdot}$ to each vertex $v_j$ until it has exactly $n$ edges. Here $\overline{v_j\cdot}$ denotes an edge with no second vertex (or an "imaginary" vertex that will not be considered for the coloring). Let the total amount of edges be given by $N$ where $n\leq N\leq Mn$. Now define $\ptM_i=e_i\in\R^M$ as the $i$-th unit vector, i.e., $\|\ptM_i-\ptM_{i'}\|_2=\sqrt{2}$ for $i\neq i'$. Furthermore, define $\chi_k\in\R^N$, $k=1,\ldots,M$ by
\begin{align*}
(\chi_k)_j=\begin{cases}
1 & \text{the j-th edge starts at }v_k, \\
-1 & \text{the j-th edge ends at }v_k,\\
0 & \text{otherwise}.
\end{cases}
\end{align*}
Note that $\chi_k$ is exactly $n$-sparse. For the $j$-th unit vector $e'_j\in\R^C$ set $\ptP_{j,k}=e'_j\otimes\chi_k\in\R^{CN}$ where $\otimes$ is the Kronecker product. (For simplicity we keep the two-dimensional indexing of $\ptP_{j,k}$ instead of reordering it into a single index.) We obtain
\begin{align}\label{eq:colDist}
\|\ptP_{j,k}-\ptP_{j',k'}\|_2=
\begin{cases}
\sqrt{2n} & j\neq j',\\
\sqrt{2n} & j=j' \text{ but }\overline{v_kv_{k'}}\not\in G,\\
\sqrt{2n+2} & j=j' \text{ and }\overline{v_kv_{k'}}\in G,
\end{cases}
\end{align}
i.e., for different colors or vertices which are not connected we obtain $\sqrt{2n}$. Now choose $\p=\infty$ and $\pp=\sqrt{n}$, then (\ref{eq:condLipp}) only holds for pairs $(j,k),(j',k')$ with either different colors $j\neq j'$ or not connected vertices $\overline{v_kv_{k'}}\not\in G$. Set $\dicM$ as identity matrix and $\resM\in\R^{CM\times M}$ to
\begin{align*}
(\resM)_{(j,k),i}=\begin{cases}
1 & k=i,\\
0 & \text{otherwise}.
\end{cases}
\end{align*}
Then, there is a feasible coloring of $G$ if and only if $\|(\dicM^*\resM)_{\iS}\|_{\lambda}=M^{1/\lambda}$ for $\lambda<\infty$ and $\iS$ solution of (\ref{eq:greedy}). By construction of $\resM$ the value $M^{1/\lambda}$ can only be achieved with $|\iS|=M$ and $(\resM)_{(j,k),i}=1$ for all $((j,k),i)\in\iS$. It follows that $(j,k)=(j,i)$ and thus the $i$-th vertex is assigned with the $j$-th color. Since $|\iS|=M$ each vertex is colored. On the other hand, each feasible coloring defines an index set $\iS\in\fS(\p,\pp)$ with $\|(\dicM^*\resM)_J\|_{\lambda}=M^{1/\lambda}$ what is maximal since $|\iS'|\leq M$ for all $\iS'\in\fS(\p,\pp)$.
\end{IEEEproof}

Following, we prove reconstruction results for GM-OMP given exact data or data obtained with noised sparsity pattern. Therefore, we need results shown in \cite{Tropp04}. We briefly summarize Theorem 3.1, Corollary 3.2 and Theorem 3.5 of this work for the reader: Given the Babel function
\begin{align*}
\mu_1(l)=\max\limits_{\Omega\subset\{1,\ldots,\nbrP\}, |\Omega|\leq l}\ \max_{\omega\not\in\Omega}\|(\dicM^*\dicM)_{\Omega,\omega}\|_1,
\end{align*}
OMP recovers the $L$-sparse solution of (\ref{eq:repProb}) in the noiseless case ($\varepsilon=0$) whenever $\mu_1(L)<1-\mu_1(L-1)$. There exists a weak version of OMP that chooses an index $i$ in each iteration such that $|(\dicM^*\resV)_i|\geq\lambda\max|\dicM^*\resV|$ with a weakness constant $\lambda\leq1$ holds. Weak OMP recovers an $L$-sparse solution whenever $\mu_1(L)<\lambda(1-\mu_1(L-1))$.

Let $\resM^l$ be the residual matrix after $l$ iterations and $\iS_l\in\fS$ the greedy choice returned by Alg. \ref{alg:GMOMPchoice}. We calculate the weakness parameter of GM-OMP by
\begin{align}\label{eq:GMOMPweakness}
\lambda=\min\left\{\frac{\left|(\dicM^*\resM^{l-1})_{i,j}\right|}{\|(\dicM^*\resM^{l-1})_{\cdot,j}\|_{\infty}}\ ,\ (i,j)\in\iS_l,\ l=1,\ldots,L\right\}.
\end{align}
Note that $\lambda$ can be very small depending on the range of amplitudes. Exemplary, consider
\begin{align*}
\dicM^*\resM^{l-1}=\begin{pmatrix}
1000 & 999 & \ldots & 1 \\
1 & 2 & \ldots & 1000 \\
\end{pmatrix}
\end{align*}
with $\p=\infty$, $\pp=0$. The greedy choice would either select the first or the second row. In both cases we obtain $\lambda=\min\{(1001-k)/k,\ k=1,\ldots,1000\}$ and hence $\lambda=1/1000$. However, $\lambda$ will be close to $1$ for feasible sets where $\min|(\solM)_{\iS_l}|>\max|(\solM)_{\iS_{l'}}|$ for $l<l'$ holds. Before we state our first theorem, we need the following definition.
\begin{defi}
For $L$ sets $\iS_1,\ldots,\iS_L\in\fS$ we say that the sets are $(\p,\pp)$-intersecting if there exists $(i,j)\in\iS_l$, $(i',j')\in\iS_{l'}$, $l\neq l'$ such that (\ref{eq:condSupp}) and (\ref{eq:condLipp}) hold.
\end{defi}
Let $\solM$ be the unique $L$-column-sparse solution of $\splM=\dicM\solM$ and $\supp\solM=\cup_{l=1}^L\iS_l$ with $\iS_l\in\fS(\p,\pp)$. We state
\begin{thm}\label{thm:exact}
GM-OMP recovers $\solM$ and all elements $\iS_1,\ldots,\iS_L$ of the sparsity pattern in $L$ iterations whenever $\iS_1,\ldots,\iS_L$ are not $(\p,\pp)$-intersecting and $\mu_1(L)<\lambda(1-\mu_1(L-1))$.
\end{thm}
\begin{IEEEproof}
The reconstruction of the support of $\solM$ follows directly by the exactness of weak OMP. Since the sets are not $(\p,\pp)$-intersecting the indices selected by Alg. \ref{alg:GMOMPchoice} belong to the same set $\iS_l$. Suppose an index in $\iS_l$ has not been selected. Then the matrix $C$ from Alg. \ref{alg:GMOMPchoice} is not zero and Alg. \ref{alg:GMOMPchoice} will not stop iterating. Thus, the complete set $\iS_l$ is recovered.
\end{IEEEproof}
Theorem \ref{thm:exact} has two disadvantages. First, the conditions depend on the solution and hence cannot be checked beforehand. Second, as we have seen $\lambda$ can be small and thus the conditions are quite restricting. To overcome the last problem, we define $\beta=\mu_1(L)/(1-\mu_1(L-1))$. Now, we can use an adaptive threshold in Alg. \ref{alg:GMOMPchoice} by selecting only indices such that $\lambda>\beta$ holds, i.e., in the $l$-th iteration only indices $(i,j)$ with
\begin{align*}
\frac{\left|(\dicM^*\resM^{l-1})_{i,j}\right|}{\|(\dicM^*\resM^{l-1})_{\cdot,j}\|_{\infty}}>\beta
\end{align*}
will be selected. We can state
\begin{thm}\label{thm:exactWeaker}
Using this threshold strategy in Alg. \ref{alg:GMOMPchoice}, GM-OMP recovers the solution $\solM$ in $L'$ iterations where $L\leq L'\leq\nbrP L$ whenever $\beta\leq1$. Furthermore, let $\iS'_1,\ldots,\iS'_{L'}$ be the feasible sets selected in each iteration and $\supp\solM=\cup_{l=1}^L\iS_l$. If $\iS_1,\ldots,\iS_L$ are not $(\p,\pp)$-intersecting, then there exists a partition $L'_1,\ldots,L'_l$ of $\{1,\ldots,L'\}$ such that $\iS_l=\cup_{l'\in L'_l}\iS'_{l'}$.
\end{thm}
\begin{IEEEproof}
Note that the first index selected in Alg. \ref{alg:GMOMPchoice} is the maximum of $C=|\dicM^*\resM|$ and thus equivalent to an OMP choice, i.e., Alg. \ref{alg:GMOMPchoice} selects at least one element per iteration. As $\beta\leq1$ this choice is part of $\supp\solM$. The algorithm terminates after at most $\nbrP L$ iterations what is the number of non-zero entries in $\solM$. If the sets $\iS_1,\ldots,\iS_L$ are not $(\p,\pp)$-intersecting the indices selected by Alg. \ref{alg:GMOMPchoice} belong to the same set $\iS_l$. Since we use a threshold, we can no longer follow that $\iS_l$ is found in one iteration. However, because $\supp\solM$ is recovered completely after $L'$ iterations, the existence of such a partition follows.
\end{IEEEproof}
Now, let us discuss the reconstruction qualities of GM-OMP due to noised sparsity patterns. Instead of $\splM=\dicM\solM$ with $\supp\solM=\cup_{l=1}^L\iS_l$, the noised data $\n\splM=\dicM\n\solM$ with $\supp\n\solM=\cup_{l=1}^L\n\iS_l$ is given. Here $\n\iS_l$ is the noised version of the pattern $\iS_l$. We consider two different kinds of noise and analyze in which case GM-OMP is able to reconstruct the sets $\n\iS_l$, $l=1,\ldots,L$. Given $\n\iS_l$ we can try to recover $\iS_l$ using a post processing strategy that will be introduced in the next subsection.
\begin{thm}\label{thm:unoise}
Let $\iS_l$ be corrupted by uniform noise $\varepsilon_u>0$:
\begin{align*}
(i,j)\in\n\iS_l && \Rightarrow && (i,j')\in\iS_l,\ d_{\domP}(\ptP_j,\ptP_{j'})\leq\varepsilon_u.
\end{align*}
Set $m=\min_{i\neq i'}d_{\domM}(\ptM_i,\ptM_{i'})$. If (\ref{eq:condLipp}) holds for $\iS_l$ with parameter $\pp$, then (\ref{eq:condLipp}) holds for $\n\iS_l$ with Lippschitz parameter $\n\pp=\pp+2\varepsilon_u/m$. If $\iS_1,\ldots,\iS_L$ are not $(\p,\pp+4\varepsilon/m)$-intersecting, then $\n\iS_1,\ldots,\n\iS_L$ are not $(\p,\pp+2\varepsilon/m)$-intersecting.
\end{thm}
\begin{IEEEproof}
For $(i_1,j_1),(i_2,j_2)\in\n\iS_l$ and $(i_1,j'_1),(i_2,j'_2)\in\iS_l$
\begin{align*}
d_{\domP}(\ptP_{j_1},\ptP_{j_2})
&\leq d_{\domP}(\ptP_{j_1},\ptP_{j'_1}) + d_{\domP}(\ptP_{j'_1},\ptP_{j'_2}) + d_{\domP}(\ptP_{j'_2},\ptP_{j_2}) \\
&\leq (\pp+2\varepsilon_u/m)d_{\domM}(\ptM_{i_1},\ptM_{i_2})
\end{align*}
holds. Equivalently, $(i_1,j_1)\in\n\iS_l$, $(i_2,j_2)\in\n\iS_{l'}$, $(i_1,j'_1)\in\iS_l$, $(i_2,j'_2)\in\iS_{l'}$ with $l\neq l'$ using that $\iS_l,\iS_{l'}$ are not $(\p,\pp+4\varepsilon/m)$-intersecting and the inverse triangle inequality:
\begin{align*}
d_{\domP}(\ptP_{j_1},\ptP_{j_2})
\geq d_{\domP}(\ptP_{j'_1},\ptP_{j'_2})-2\varepsilon_u
>(\pp+\frac{2\varepsilon}{m})d_{\domM}(\ptM_{i_1},\ptM_{i_2}).
\end{align*}
\end{IEEEproof}
The noise assumed in Theorem \ref{thm:unoise} typically appears in applications where measurements may be corrupted due to shaking apertures. If an upper bound $\varepsilon_u$ is known, the parameters of GM-OMP can be adapted.
\begin{thm}\label{thm:bnoise}
Let $\iS_l$ be corrupted by Bernoulli distributed noise $\varepsilon_B\in[0,1]$, i.e.,
\begin{align*}
\n\iS_l\subseteq\iS_l, && \mathop{Pr}\left((i,j)\not\in\n\iS_l\ |\ (i,j)\in\iS_l\right)=\varepsilon_B
\end{align*}
where $\mathop{Pr}\left((i,j)\not\in\n\iS_l\ |\ (i,j)\in\iS_l\right)$ is the probability that an index $(i,j)\in\iS_l$ is not in the corrupted set $(i,j)\not\in\n\iS_l$. Let (\ref{eq:condSupp}) hold for $\iS_l$ with parameter $\p$ and $\mathop{Pr}(\n\iS_l\in\fS(k\p,\pp))$ be the probability that (\ref{eq:condSupp}) holds for $\n\iS_l$ with parameter $k\p$, $k\in\N$. Then
\begin{align*}
(1-\varepsilon_B^k)^{\nbrM-k+1}\leq\mathop{Pr}(\n\iS_l\in\fS(k\p,\pp)).
\end{align*}
\end{thm}
\begin{IEEEproof}
Note that (\ref{eq:condSupp}) gives a connected graph. We search for a lower bound of the probability, that the graph is still connected when we remove points $\ptM_j$ with probability $\varepsilon_B$ but add edges $\overline{\ptM_j\ptM_{j'}}$ with $d_{\domM}(\ptM_j,\ptM_{j'})\leq k\p$. For a lower bound it is sufficient to consider the worst case, i.e., $\ptM_j=j\p$. The graph becomes a line with at most $\nbrM$ points. The graph is connected whenever there is a connection from $\ptM_1$ to $\ptM_{\nbrM}$. For the new parameter $k\p$ we obtain the edges $\overline{\ptM_j,\ptM_{j+q}}$ with $q\leq k$. It follows that the graph will no longer be connected whenever $k$ consecutive points vanish.

This problem is an application of success runs in Bernoulli trails \cite{Muselli00}. In particular, $\mathop{Pr}(\n\iS_l\in\fS(k\p,\pp))$ is the probability that the longest success run is shorter than $k$. This probability has an exact but rather complicated analytic expression. The simple lower bound $(1-\varepsilon_B^k)^{\nbrM-k+1}$ is shown in \cite{Fu85}. Other bounds and the exact analytic form can also be found in \cite{Muselli00}.
\end{IEEEproof}
The noise assumed in Theorem \ref{thm:bnoise} appears in applications e.g., whenever a single measurement is lost or a sensor fails. The parameter $\p$ can be adapted according to Theorem \ref{thm:bnoise}.

Theorem \ref{thm:exact} gives two conditions for exact recovery. The condition $\mu_1(L)<\lambda(1-\mu_1(L-1))$ ensures recovery of the right support set. It was deduced from OMP and was shown to be strict \cite{Tropp04}. The $(\p,\pp)$-separation condition guarantees the separation of the support into $L$ structures. This condition is not strict. The $L$ feasible sets $\iS_1,\ldots,\iS_L$ may still be reconstructed without having $(\p,\pp)$-separation depending on the amplitudes $(\solM)_{\iS_l}$, $l=1,\ldots,L$. Theorem \ref{thm:exactWeaker} gives reconstruction results if one or both of these conditions were penalized.

In Theorem \ref{thm:unoise} and \ref{thm:bnoise} we discussed a noised sparsity pattern and how the parameters $\p,\pp$ should be adapted. In the next section, we present a post processing step to reconstruct the original pattern given. Beforehand, we give a statement on two other cases of noisy data. First, consider the case where $\splM$ does not have an exact sparse representation $\solM$ with $\dicM\solM=\splM$, but instead we search for a sparse approximation as in problem \ref{eq:repProbMMV}. This problem occurs e.g., when the data $\splM$ is noised. We can easily obtain similar results to Theorem \ref{thm:exact} and \ref{thm:exactWeaker} by replacing the exact recovery condition of (weak) OMP with the optimal $L$-term approximation conditions given in \cite{Tropp04}.

As another scenario, consider a sampling $\splM$ that is sparse in some dictionary $\dicM$, i.e., there exists a sparse solution $\solM$ of $\dicM\solM=\splM$. Now, instead of $\dicM$ we have only given the dictionary $\n\dicM$. Exemplary, let $\splM$ be sparse in Fourier domain but not necessarily containing frequencies given by the discrete Fourier transform. Given the Fourier transform of $\splM$ is it possible to reconstruct the exact frequencies, i.e., given an approximation in $\n\dicM$ is it possible to reconstruct the exact dictionary $\dicM$? This problem was analyzed under the keyword of super-resolution in \cite{Candes14} for the SMV problem. Only recently, the MMV problem with common support constraint was discussed in \cite{Yang16}. In both cases, an exact recovery is possible whenever the non-zero entries are separated by at least a distance depending on the super-resolution factor. An interesting question we consider for future work, is the connection between this separation and patterns that are not $(\p,\pp)$-intersecting. This connection may be used to design a super-resolution method for generalized patterns.

\subsection{Post processing}

So far, we presented the GM-OMP algorithm and proved basic theoretical properties. Before we demonstrate the technique on numerical examples, we discuss how to use GM-OMP for powerful post processing of the data. Consider we reconstructed a solution $\solM$ and its support $\supp\solM=\setI=\n\iS_1\cup\ldots\cup\n\iS_L$, where $\n\iS_l\in\fS$ is assumed to be a corrupted sparsity pattern.

While it is a common idea to denoise corrupted amplitude values of $\solM$, the sparsity pattern has been of minor interest so far. Even though the pattern itself might be noised. Exemplary, in non-destructive testing external forces during the measurement can corrupt the probes positions what influences the geometry and thus the sparsity pattern \cite{Kaaresen99,Bossmann12}. As another example, consider $\dicM$ being the Fourier matrix. It only contains a fixed amount of Fourier frequencies. However, there are signals that are sparse in Fourier domain but only consist of frequencies not covered by the matrix. Then the reconstructed sparse approximation most likely rounds these frequencies upto the closest frequency of $\dicM$, what can be interpreted as a corrupted sparsity pattern of $\solM$. As last example, simply assume a failed measurement, i.e., a zero column in $\splM$. Surely the corresponding column in $\solM$ will also be zero. To reconstruct the original signal, we can apply inpainting ideas on the sparsity pattern.

Remembering Fig. \ref{fig:suppM}, i.e., $\iS_l=\supp\mtx_l$ as a discrete sampling of a function, we can denoise the sparsity pattern for $l=1,\ldots,L$ by solving the problems
\begin{align}\label{eq:patternNoise}
f_l=\argmin\limits_{f\in\F}\left\|\left(\ptP_j-f(\ptM_i)\right)_{(i,j)\in\n\iS_l}\right\|_2+\delta|\supp f|,
\end{align}
with a weight $\delta>0$. Here, $\F=\F(\domM,\domP)$ is a suitable function space (e.g., polynomials, splines, $\ldots$). Afterwards set the denoised pattern $\iS_l$ to
\begin{align*}
\iS_l=\{(i,j)\ |\ \ptM_i\in\supp f_l,\ \ptP_j=[f_l(\ptM_i)]\}
\end{align*}
where $[f(\ptM_i)]$ is $f(\ptM_i)\in\domP$ rounded to the closest of the elements $\ptP_1,\ldots,\ptP_{\nbrP}$. For small $\delta$ the support $\supp f$ may be large and hence $|\iS_l|$ can increase. This gives an inpainting strategy to reconstruct missing structure elements.

A similar approach can be applied to denoise the amplitudes of $\solM$. Given $\iS_l$ and assume that $\iS_l\cap\iS_{l'}=\emptyset$ for all $l'\neq l$ we solve
\begin{align}\label{eq:ampNoise}
g_l=\argmin\limits_{g\in\G}\left\|\left((\solM)_{i,j}-g(\ptM_i)\right)_{(i,j)\in\iS_l}\right\|_2
\end{align}
on a function space $\G=\G(\domM,\C)$ and update $(X)_{i,j}=g_l(\ptM_i)$, for all $(i,j)\in\iS_l$.

\section{Numerics}

We demonstrate the advantages of our proposed algorithm in three examples. First, we compare the technique with other sparse approximation methods for the MMV problem. Afterwards we discuss two practical examples and illustrate the information given by the sparsity pattern.

\subsection{Numerical comparison}

\begin{figure}
\centering
\includegraphics[width=43mm]{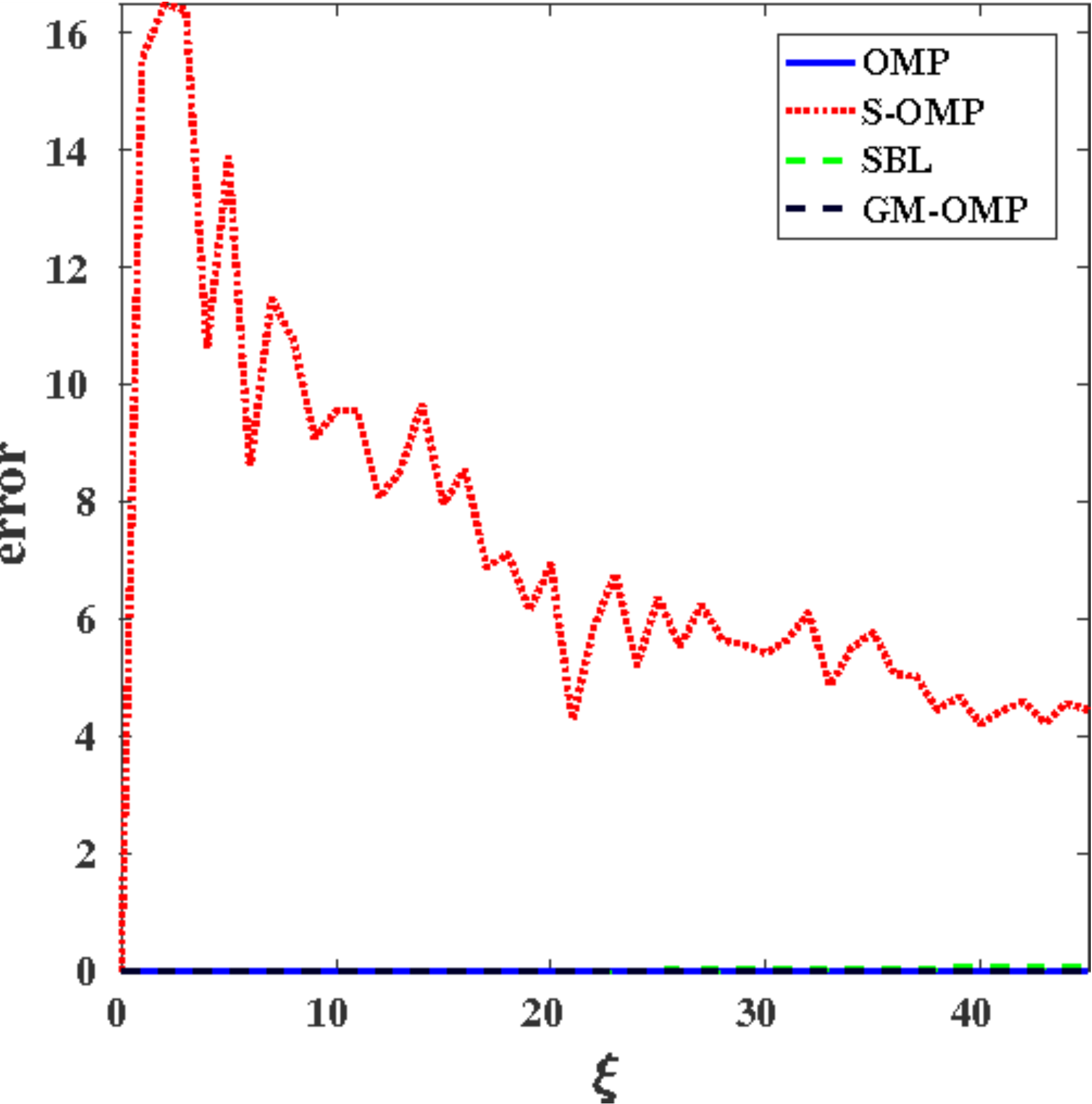}
\includegraphics[width=43mm]{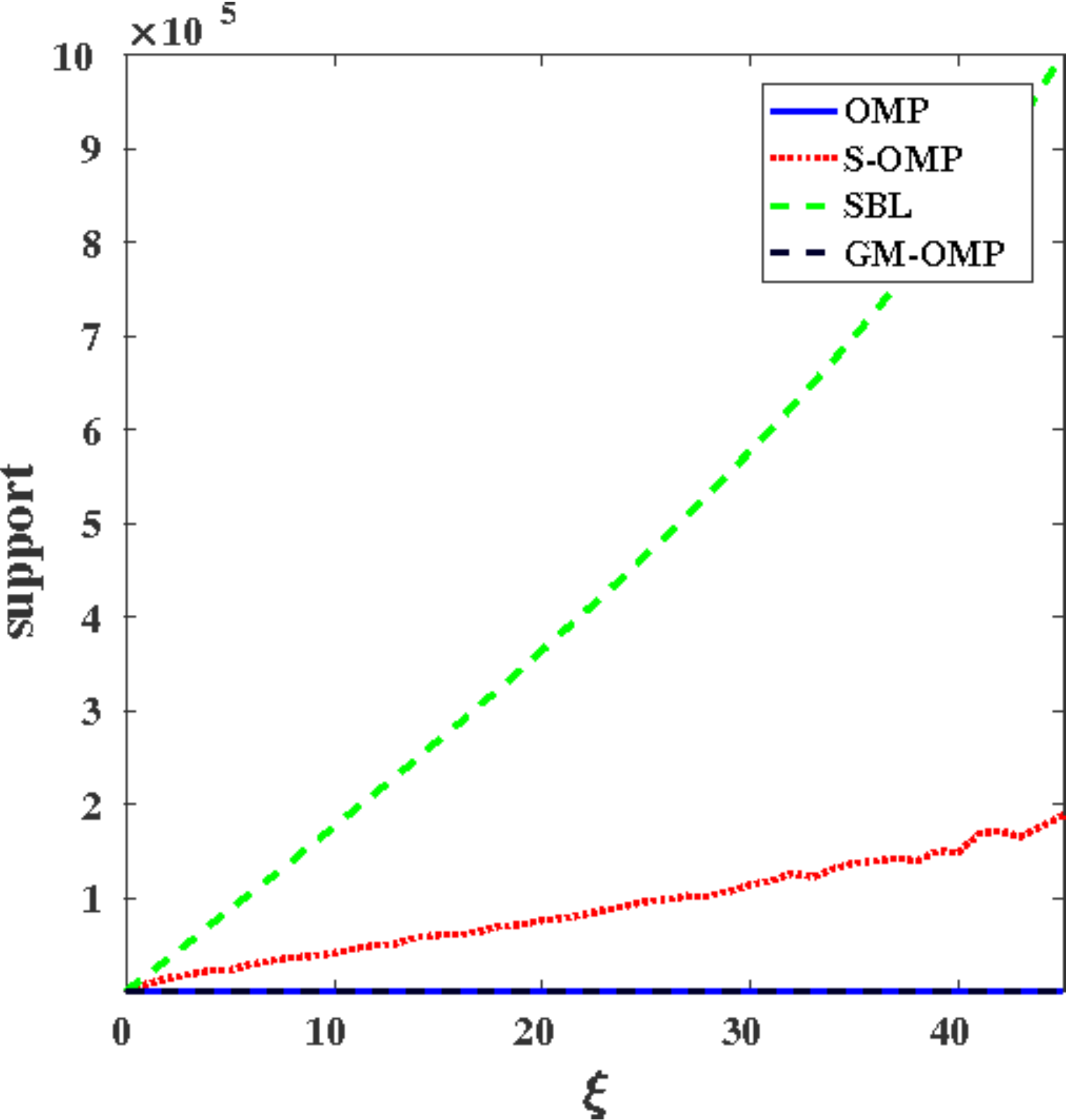}
\caption{Reconstruction error and non-zero elements of the solution.}
\label{fig:compExact}
\end{figure}

We compare our method to three other techniques: OMP applied to each column separately, S-OMP \cite{Tropp06a} and MSBL, a technique presented in \cite{Wipf07} based on sparse bayesian learning. Let $\dicM$ be a convolution matrix of a Gauss kernel with standard deviation $\sqrt{2.5}$. We define the matrix $\solM\in\R^{1000\times1000}$ by
\begin{align*}
(\solM)_{i,j}=\begin{cases}
1 & i=[j\tan\xi]\\
0 & \text{otherwise}
\end{cases}
\end{align*}
for $i,j\leq1000$, i.e., each column of $\solM$ is $1$-sparse. The matrix is clearly structured, it consists of one line  with a slope of $\xi$. For $\xi=45^{\circ}$ this becomes the identity matrix; $\xi=0^{\circ}$ gives a matrix with one non-zero row (which is the pattern that S-OMP and the MSBL assume). For $\xi=0^{\circ},\ldots,45^{\circ}$ we calculate $\splM=\dicM\solM$ and use all methods to reconstruct $\solM$. Therefore, we choose $\ptM_i=i$, $\ptP_j=j$, $\p=1$ and $\pp=1$ (which corresponds to a maximal slope of $45^{\circ}$). We use $L=1$ iteration since $\solM$ contains exactly one structure. In Fig. \ref{fig:compExact} the reconstruction error and the number of non-zero elements in the solution $\solM$ are shown for all algorithms. Nearly all methods are able to find a good approximation. Only S-OMP forces row sparsity of $\solM$ and thus produces stare casing effects which corrupt the solution. Both the MSBL and S-OMP assume that $\solM$ is row-sparse, i.e., there are only a few non-zero rows. Once a row contains a non-zero element, the entire row is considered to be non-zero. This leads to an extreme overestimation of the support while OMP and GM-OMP can find the exact number of non-zero entries.

\begin{figure}
\centering
\includegraphics[width=43mm]{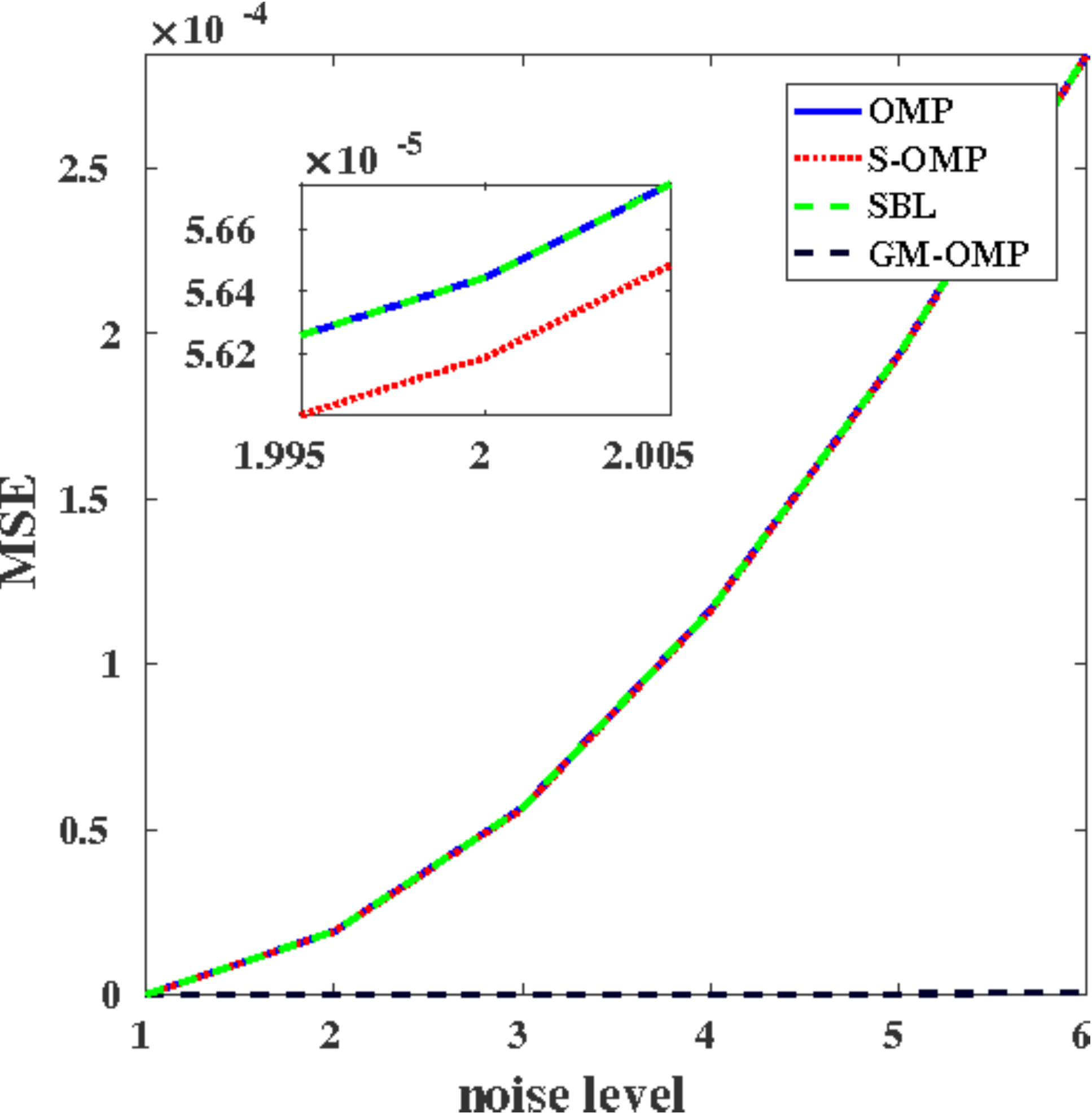}
\includegraphics[width=43mm]{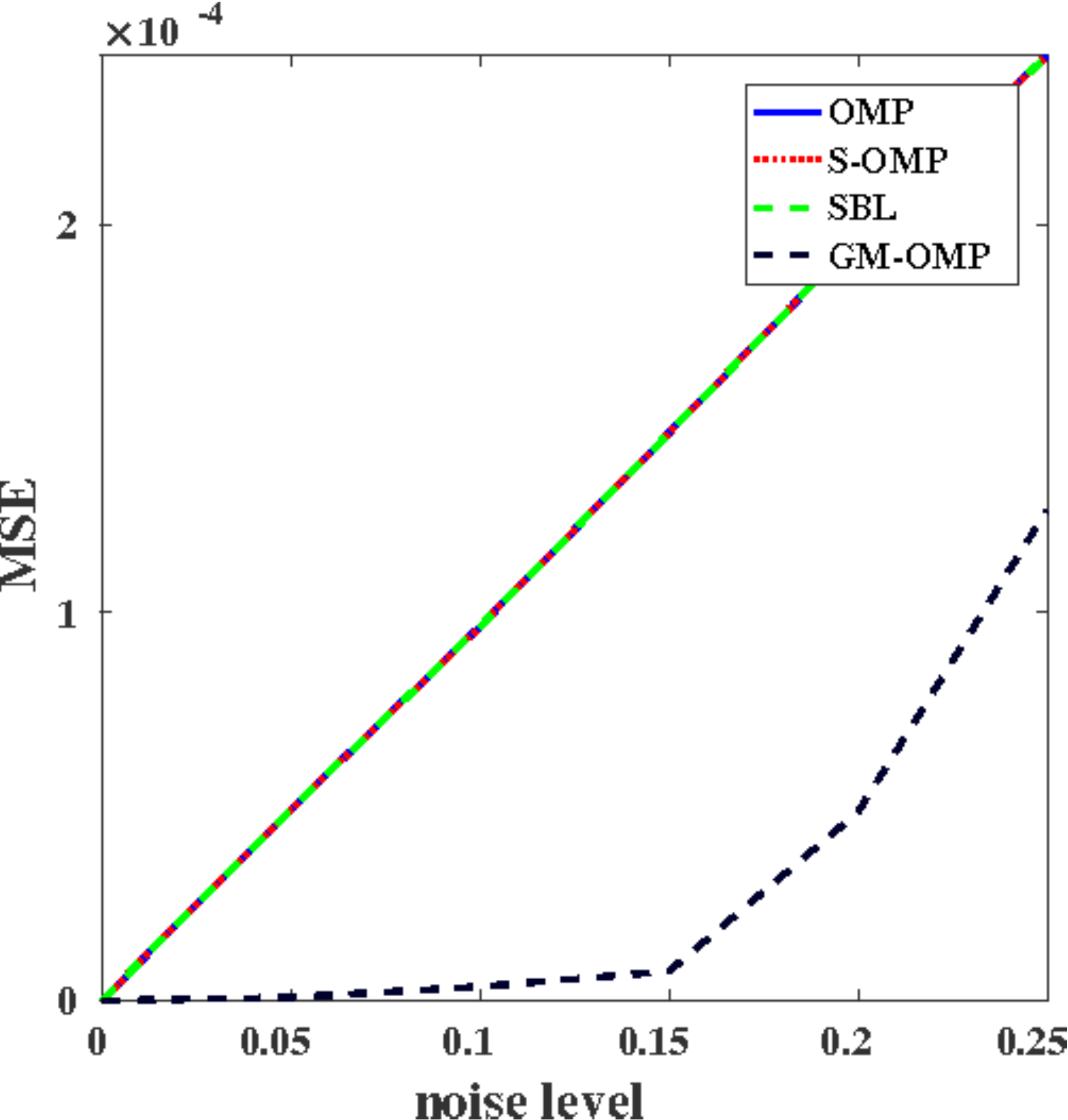}
\caption{Reconstruction error for the noisy data $\solM_1$ (left) and $\solM_2$ (right).}
\label{fig:compNoise}
\end{figure}

Next, we demonstrate the power of the proposed denoising step. Therefore consider $\solM$ and its noised versions $\solM_1,\solM_2\in\R^{1000\times1000}$ with
\begin{align*}
(\solM)_{i,j}=\begin{cases}
1 & i=500\\
0 & \text{otherwise}
\end{cases},
&&
(\solM_1)_{i,j}=
\begin{cases}
1 & i=500+[\varepsilon_u(j)]\\
0 & \text{otherwise}
\end{cases}.
\end{align*}
\begin{align*}
(\solM_2)_{i,j}=\begin{cases}
\varepsilon_B(j) & i=500\\
0 & \text{otherwise}
\end{cases},
\end{align*}

where $|\varepsilon_u(j)|\leq\varepsilon_u$ is uniform noise and $\varepsilon_B(j)\in\{0,1\}$, $\mathop{Pr}(\varepsilon_B(j)=0)=\varepsilon_B$ is Bernoulli distributed. Given $\splM_1=\dicM\solM_1$ or $\splM_2=\dicM\solM_2$ we want to reconstruct $\solM$. The mean squared error over $100$ runs for both cases is plotted in Fig. \ref{fig:compNoise}. The noise level gives the values of $\varepsilon_u$ respectively $\varepsilon_B$.

For $\solM_1$ we adapted the parameters according to Theorem \ref{thm:unoise} and set $\pp=13\geq1+2\varepsilon_u/m$. In the second case, we choose $\p=6$ for the reconstruction of $\solM_2$. Using Theorem \ref{thm:bnoise} this gives us a reconstruction probability of more than $78\%$ even for $\varepsilon_B=0.25$. In both cases we solve optimization problem (\ref{eq:patternNoise}) with $\delta=0$ and $\F=\Pi_4$ afterwards to denoise the pattern.

Note that the row-sparsity assumption of S-OMP and MSBL hold for $\solM$. Nevertheless, only GM-OMP is able to find a good approximation in most cases. All other methods only reconstruct the noised matrices $\solM_1,\solM_2$. For a high probability of $\sigma=0$ the MSE of GM-OMP increases, i.e., the parameter $\p$ does no longer compensate the missing data. Interestingly, S-OMP profits from its stare casing effects when the pattern is uniformly noised and hence returns a slightly better solution.

\subsection{Application 1: Non-destructive testing}

\begin{figure}
\centering
\includegraphics[width=43mm]{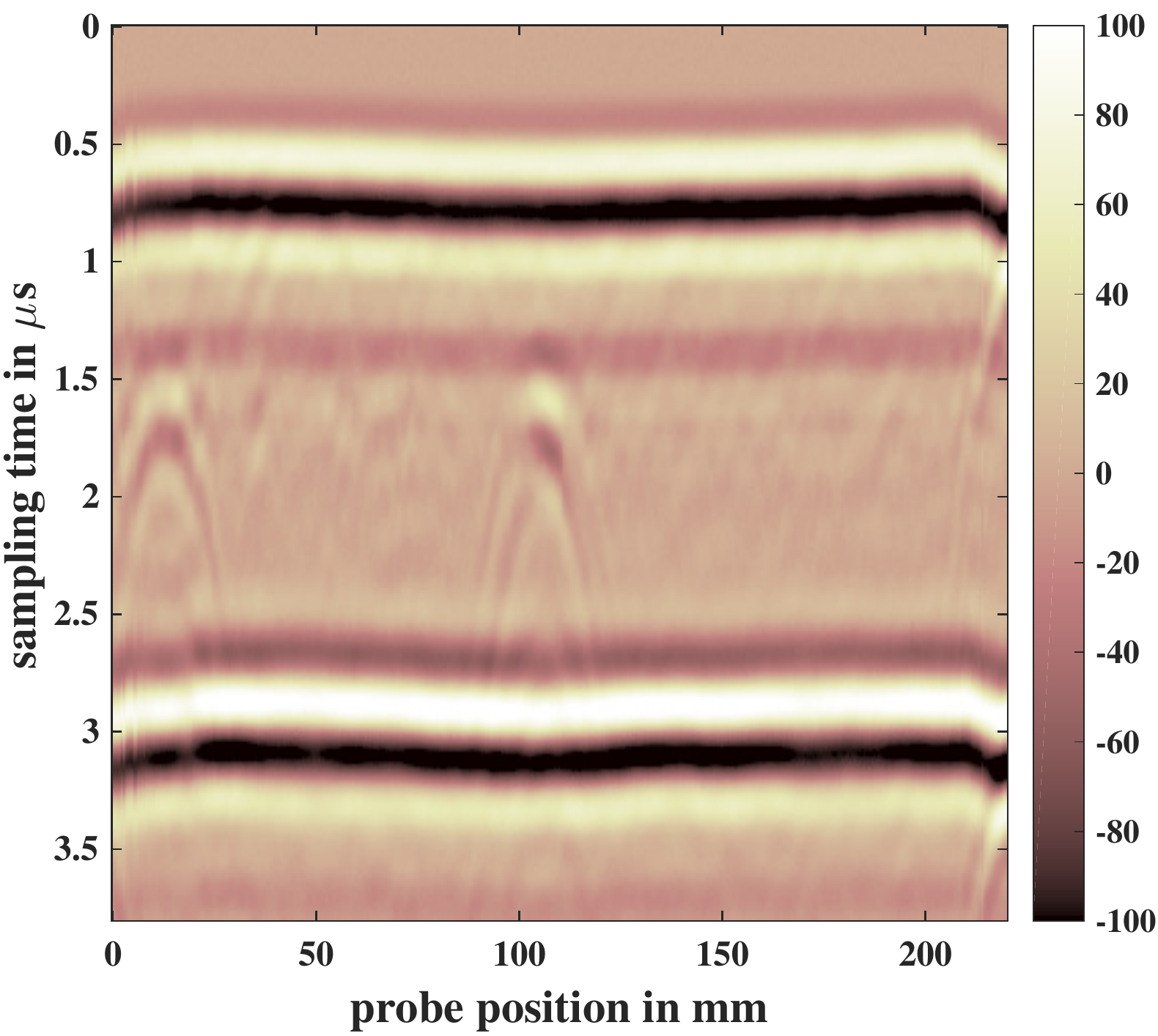}
\caption{Real ultrasonic non-destructive testing data.}
\label{fig:tofdOrg}
\end{figure}

As a first practical example we analyze ultrasonic data obtained from non-destructive testing of steel tubes. The original data shown in Fig. \ref{fig:tofdOrg} was generated by the "Time-of-Flight Diffraction" (ToFD) method using an Olympus Omniscan iX system with $5$Mhz transducer, $6$mm diameter and $70^{\circ}$ angle of incidence. The tested tube was a large diameter pipe with outer diameter $1066$mm and $23.3$mm wall thickness. Each column of the data represents a measured signal at different positions on the tubes surface. The positions were equidistantly set on a straight line with a distance of $0.5$mm, hence we set $\ptM_i=0.5i$. The signals were measured in time with a sampling ratio of $0.01\mu$s.
Four major events can clearly been seen in the data. The topmost one is an ultrasonic impulse that directly travels through the surface from transducer to receiver - the lateral wave. The bottommost is an impulse that was reflected by the back wall - the back wall echo. The two events in between (recognizable as parabolas) indicate defects in the material. We use GM-OMP to recover and denoise these events.

Ultrasonic data is column-wise sparse when $\dicM$ is a convolution matrix based on the Gabor impulse (\cite{Bossmann12,Bossmann15})
\begin{align}
g(t) = e^{-\theta t^2}\cos(\phi t+\psi).\label{eq:Gabor}
\end{align}
Here $\theta$ is the bandwidth factor, $\phi$ is the frequency and $\psi$ is the phase. Thus, $\ptP_j=0.01j$ is the shift of the $j$-th column in $\mu$s. For the given data we choose $p=6.8486$, $\phi=14.685$ and $\psi=-2.0836$ (see \cite{Bossmann12,Bossmann15} for details about the parameter choice). We define our feasible set $\fS$ using $\p=0.5$ such that$|\ptM_i-\ptM_{i'}|\leq\p$ only for $i'=i-1,i+1$. Note that (\ref{eq:condLipp}) compares distance in time ($\mu$s) with a distance in space (mm). The ultrasonic speed in steel is about $5.9$mm/$\mu$s, hence we chose $\pp=0.1>5.9^{-1}$ what gives stability for noisy data.

After applying $L=4$ iterations of GM-OMP to the data, we use the denoising strategies discussed in the last section and set $\beta=1$. Since structures in pipe testing often behave linearly or quadratically, we use $\F=\{f\cdot\chi_C\ |\ f\in\Pi_4,C\subseteq\domM\text{ convex }\}$, i.e., polynomials upto degree four multiplied by a characteristic function. The characteristic function is needed for the support constraint in (\ref{eq:patternNoise}). Moreover we choose $\G=\Pi_0$ as the space of all constant functions, i.e., the amplitudes of each structure are set to its mean value. This value can give a first hint about the underlying material in the pipe.

\begin{figure}
\centering
\includegraphics[width=43mm]{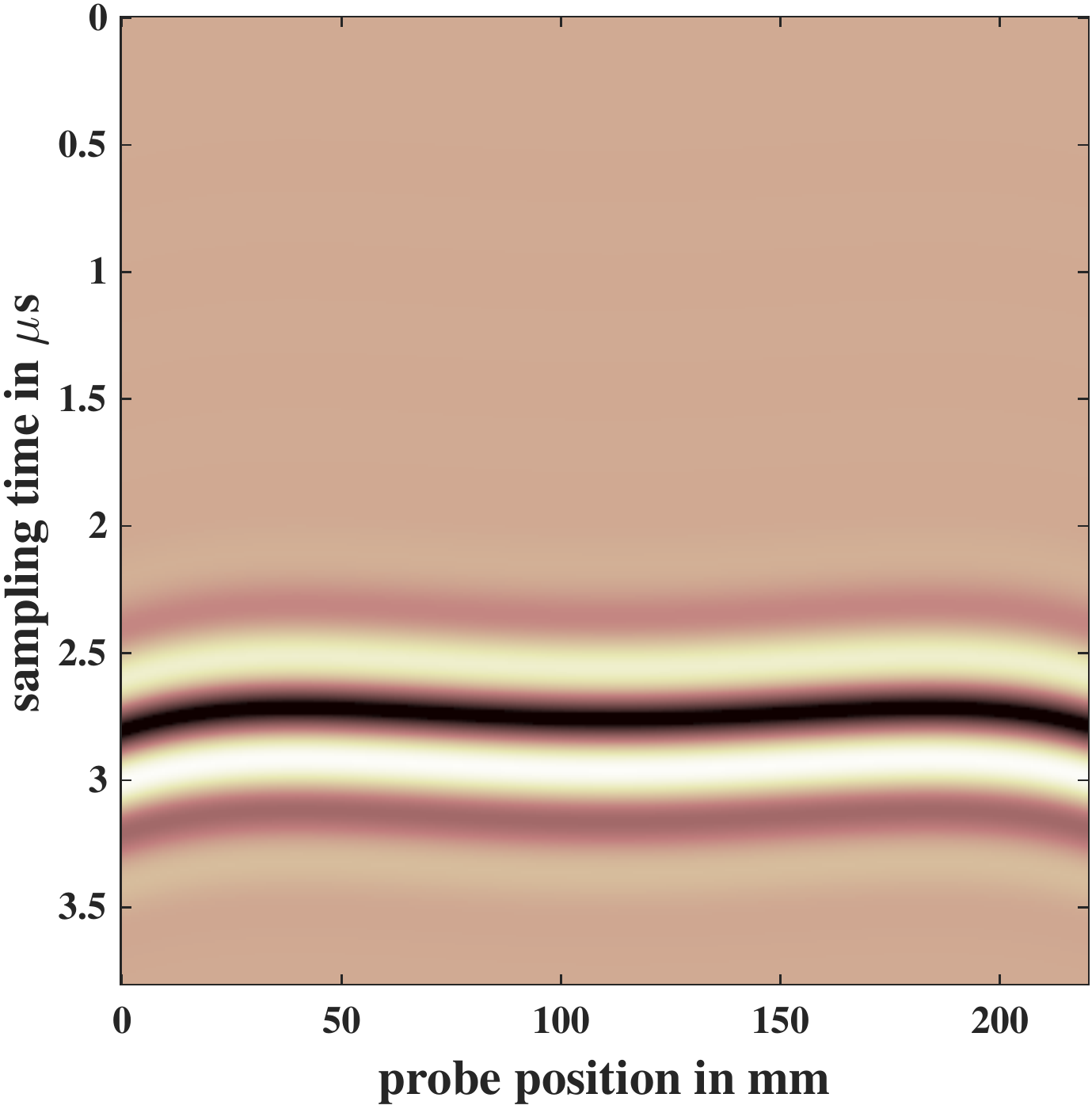}
\includegraphics[width=43mm]{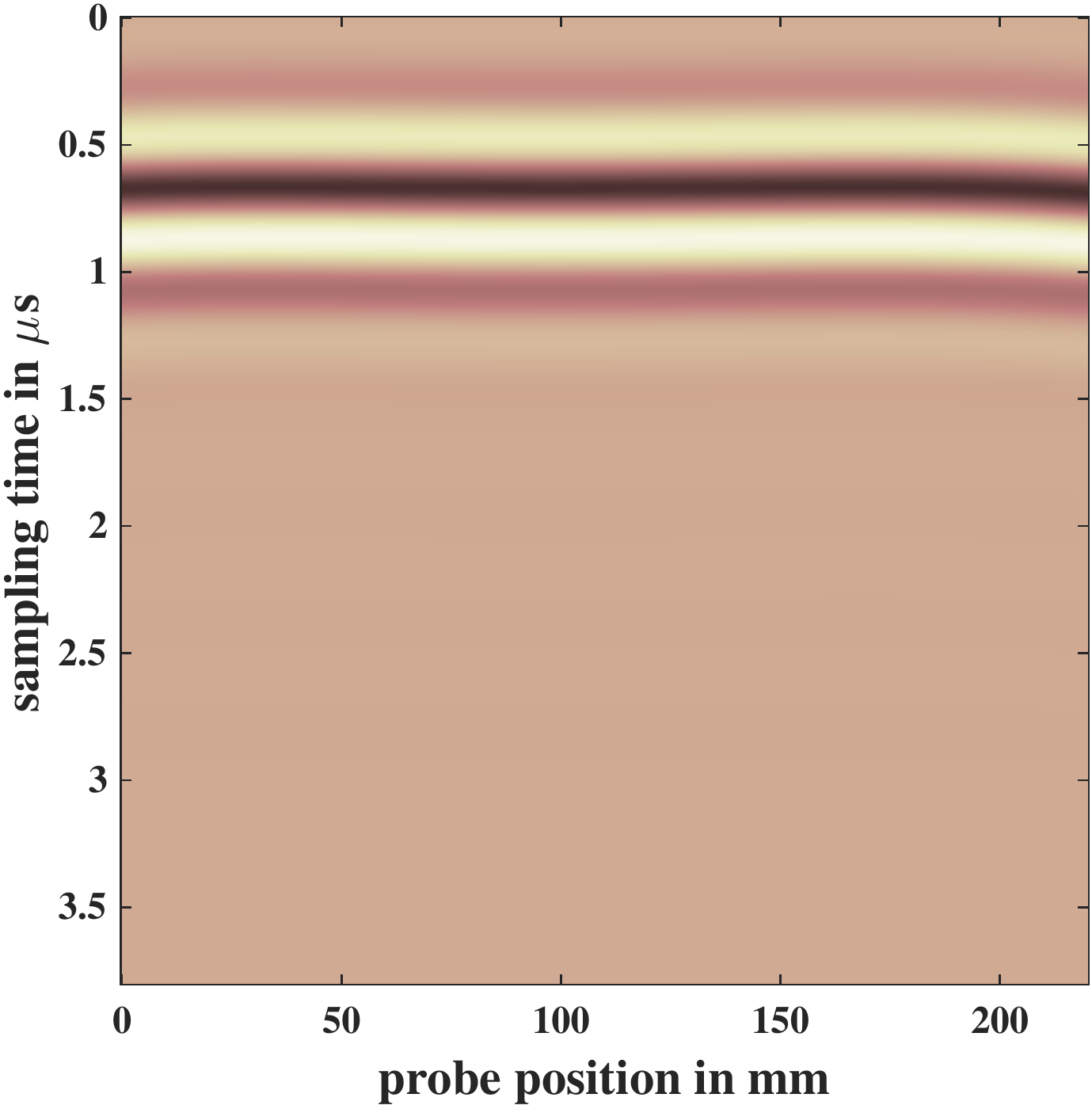}\\
\includegraphics[width=43mm]{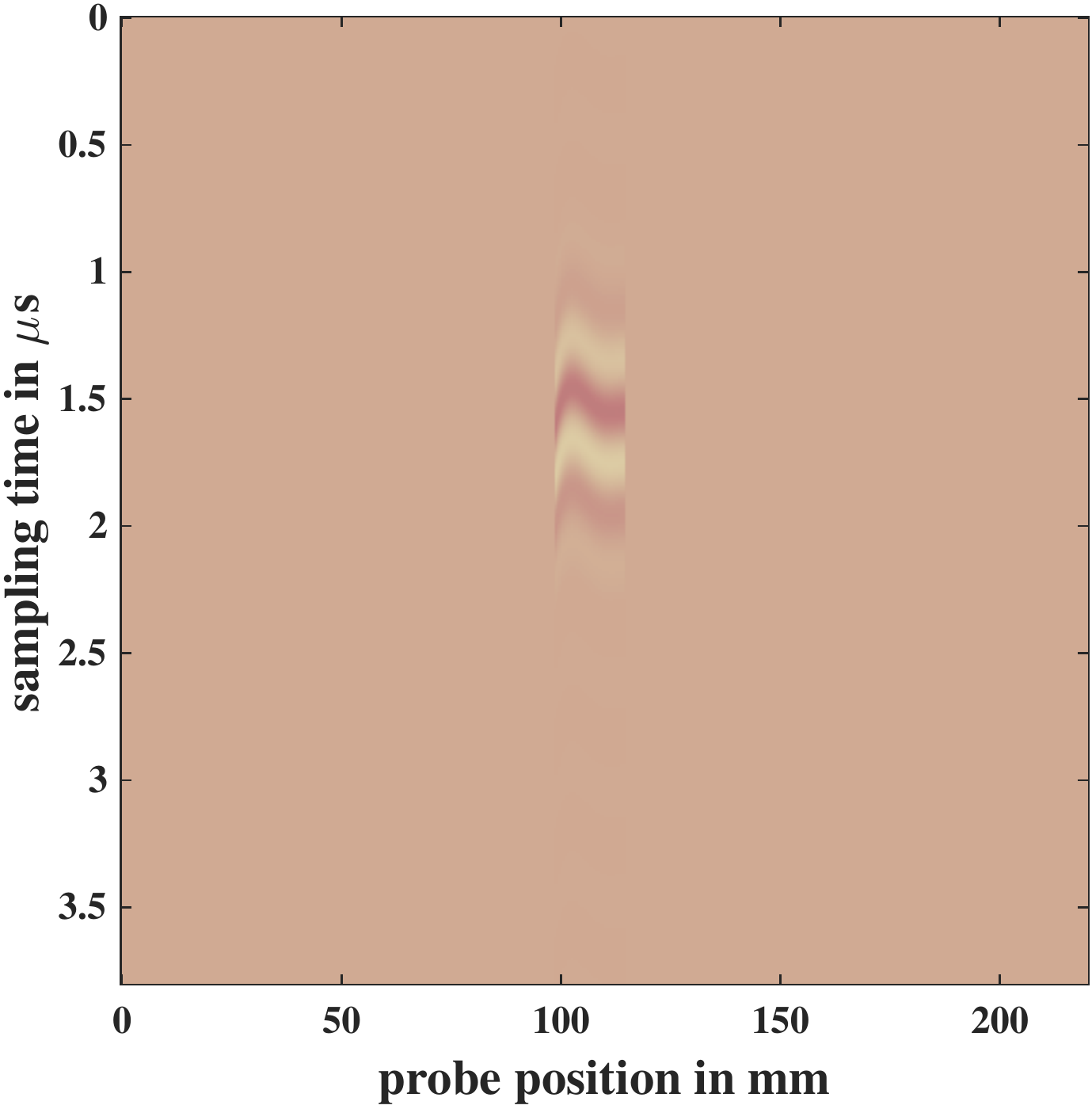}
\includegraphics[width=43mm]{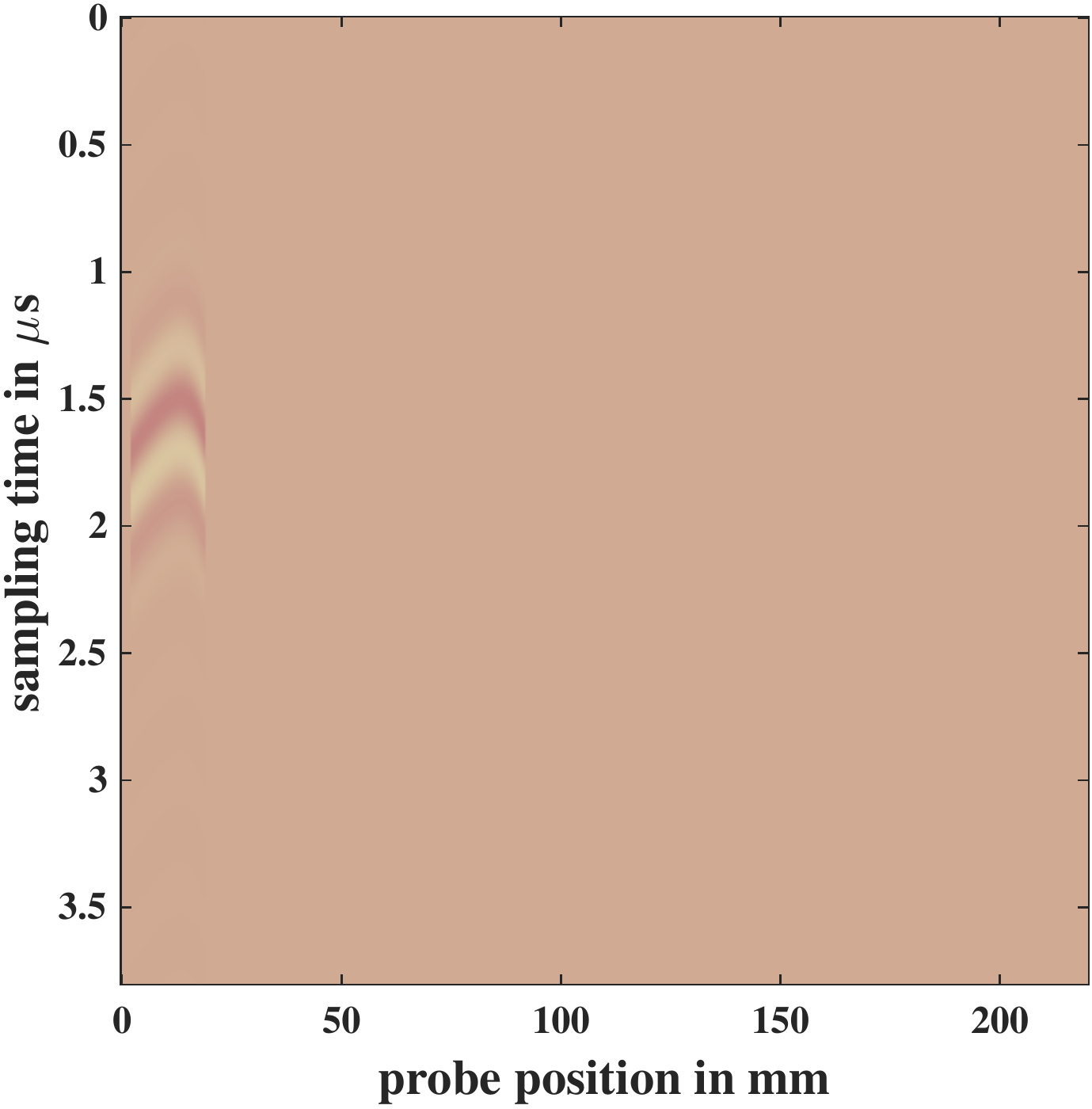}
\caption{Reconstructed and denoised structures.}
\label{fig:tofdE}
\end{figure}

In Fig. \ref{fig:tofdE} the four sparsity patterns found by GM-OMP are shown in data domain (i.e., multiplied by $\dicM$). As we see, the algorithm is able to reconstruct all four structures of the original data. Due to the structural denoising, the pattern looks more smooth and effects caused by shaking probes are no longer visible. In Tab. \ref{tab:coeff} the polynomial coefficients of each pattern are shown. As suppoesd, the lateral wave and back wall echo are mostly linear while the defects were approximated by a quadratic polynomial.

\begin{table}
\caption{Polynomial coefficients of the structures shown in Fig. \ref{fig:tofdE}.}
\centering
\begin{tabular}{|c|r|r|r|r|r|}
\hline
$L$&$1$&$2$&$3$&$4$\\
\hline
$x^4$&$1.05e-10$&$1.57e-11$&$1.48e-6$&$5.65e-7$\\
$x^3$&$8.91e-08$&$1.45e-08$&$-3.02e-5$&$4.22e-5$\\
$x^2$&$2.42e-05$&$4.80e-06$&$-0.000290$&$0.001354$\\
$x$&$0.002202$&$0.000576$&$0.007798$&$0.010080$\\
const.&$2.8598$&$0.770270$&$1.6013$&$1.6006$\\
\hline
\end{tabular}
\label{tab:coeff}
\end{table}

\subsection{Application 2: meteorologic data}

\begin{figure}
\centering
\includegraphics[width=43mm]{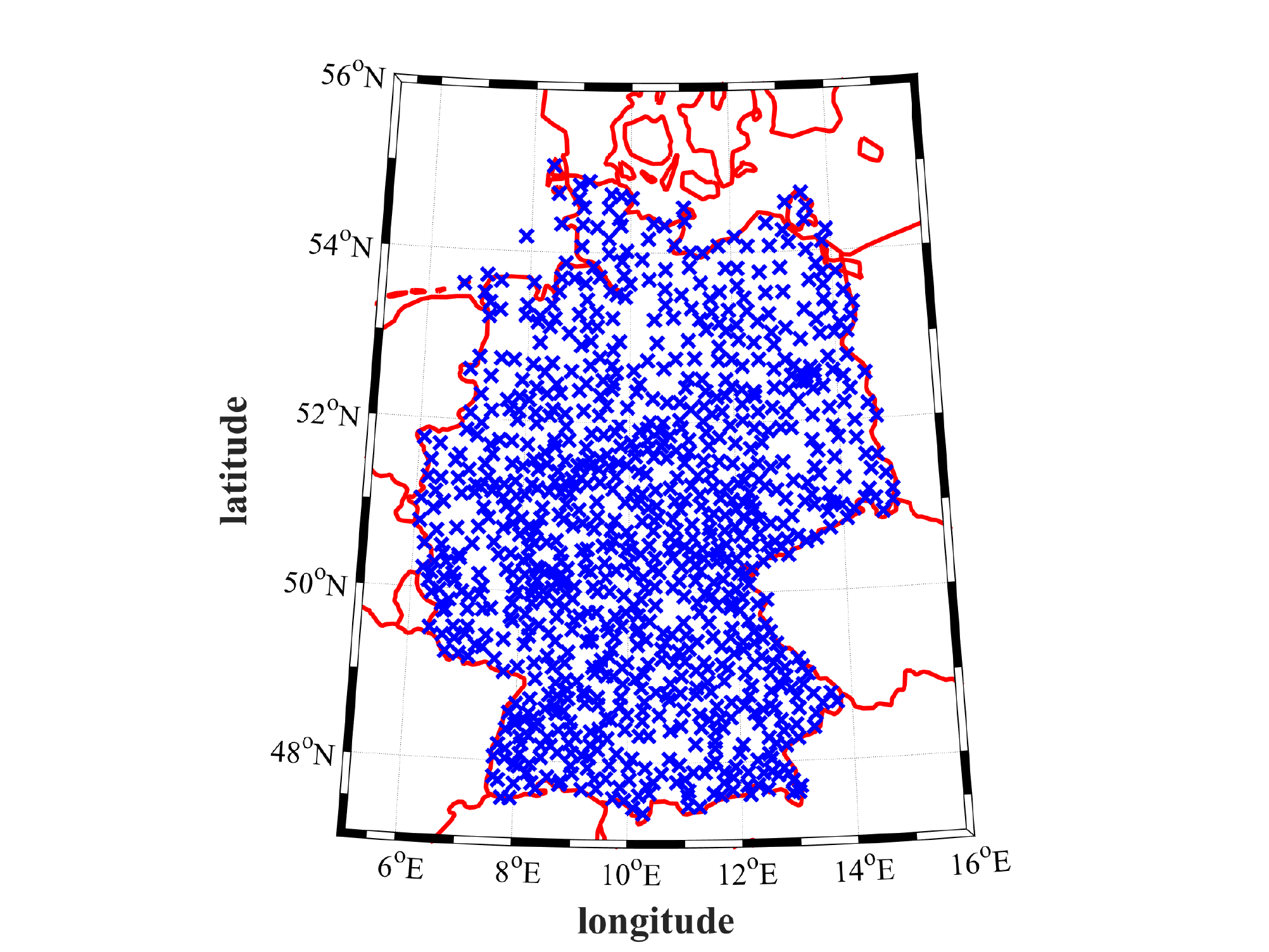}
\caption{Weather stations for precipitation in Germany \cite{MMAP}.}
\label{fig:stationMap}
\end{figure}

\begin{figure}
\centering
\includegraphics[width=43mm]{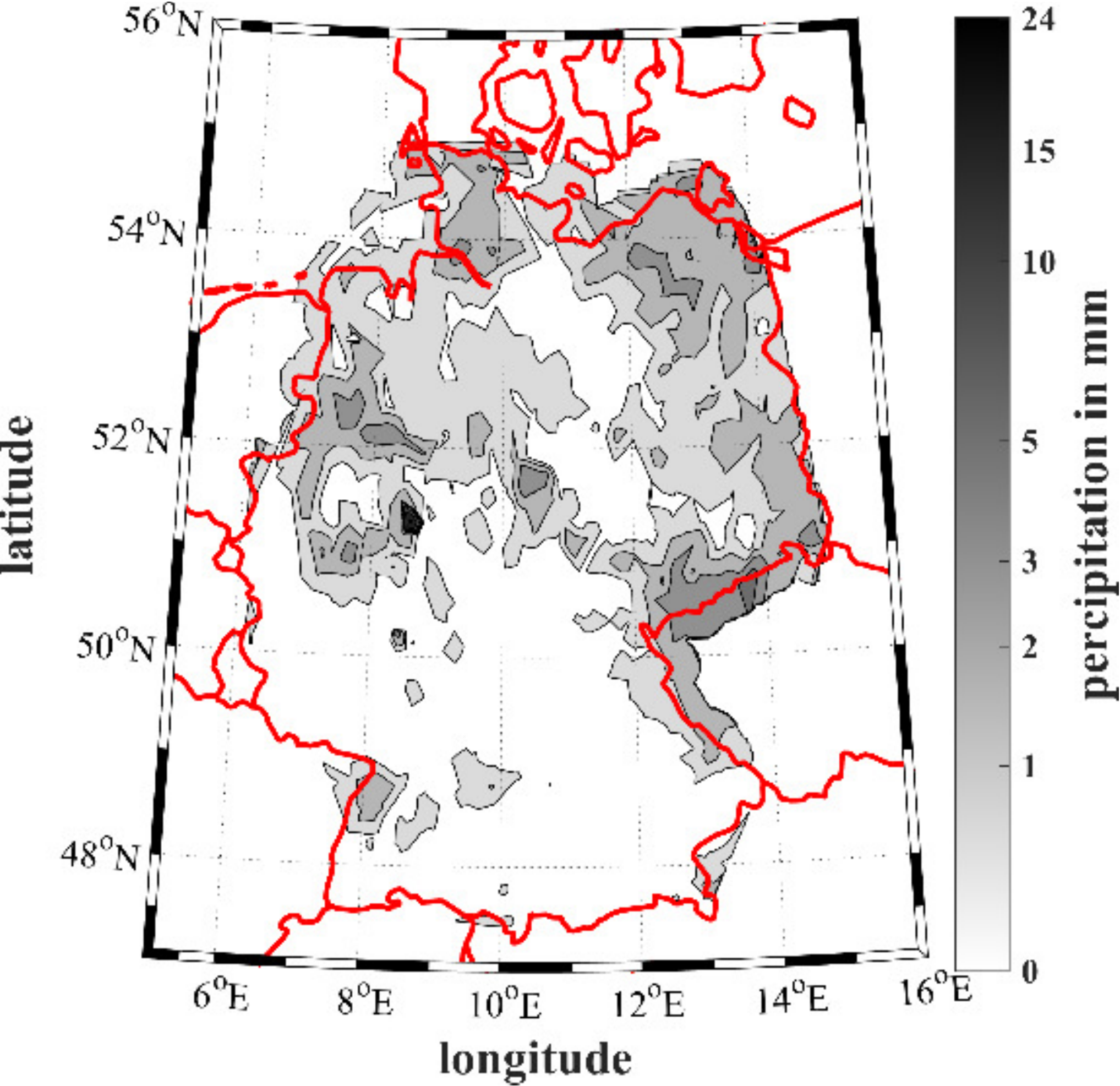}
\includegraphics[width=43mm]{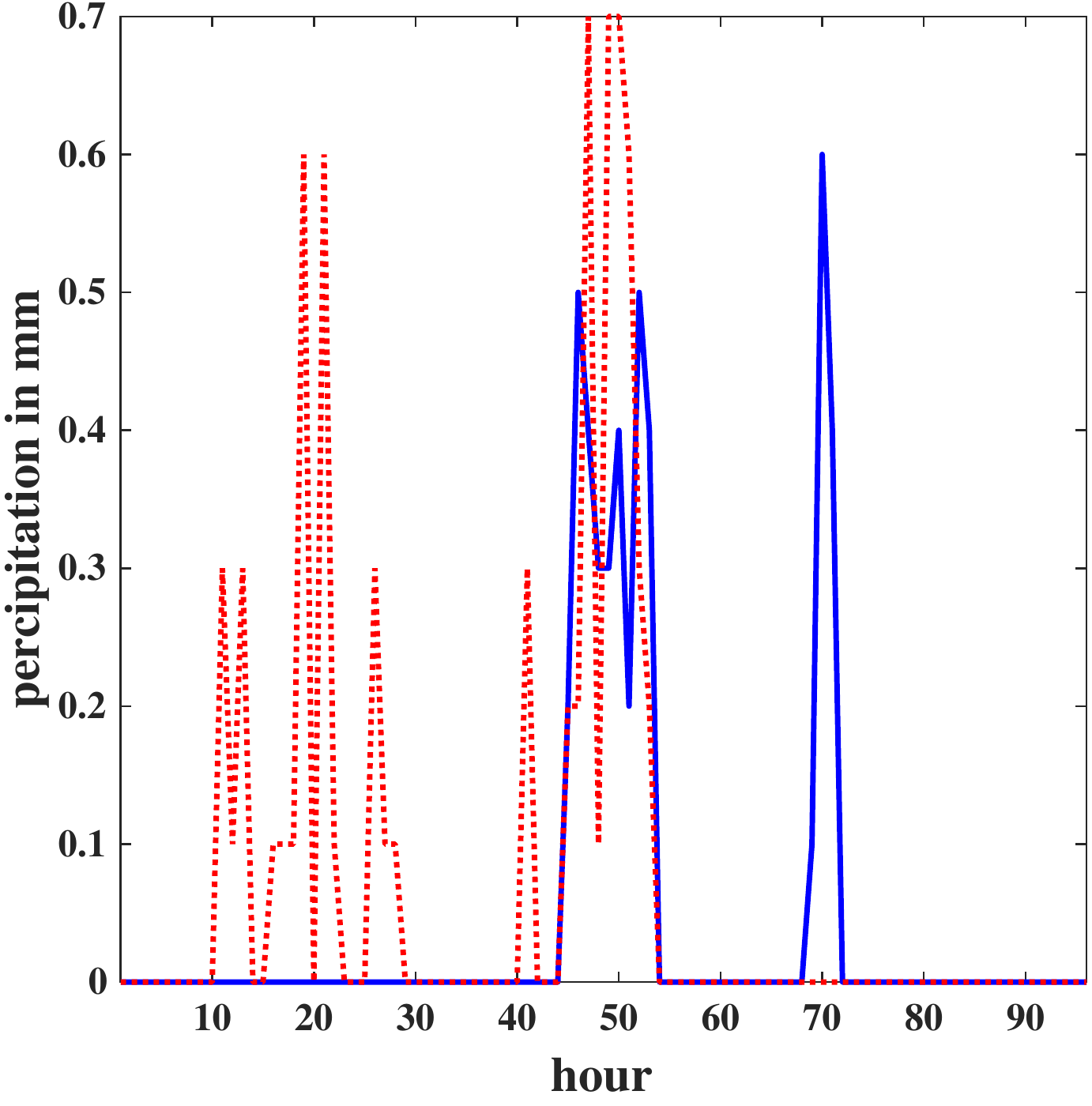}
\caption{Left: Overall precipitation (log scale) \cite{MMAP}; right: data of two stations exemplary.}
\label{fig:metData}
\end{figure}

In this example we use GM-OMP without post-processing on meteorologic data provided by Deutscher Wetterdienst (DWD) \cite{DWD}. We analyze the hourly precipitation in Germany from 25th to 28th of November 2008 where we use data of $932$ stations shown in Fig. \ref{fig:stationMap}. Stations that were moved during this time or had too many missing values were neglected. In Fig. \ref{fig:metData} the overall precipitation and data of two stations is plotted exemplary where $0$ hour refers to Nov. 25th 2008, $0:00$. We use a dictionary based on centered cardinal B-splines
\begin{align*}
B_1(t)=\chi_{[-0.5,0.5)} && B_n(t)=\int\limits_{-0.5}^{0.5}B_{n-1}(t-\tau)d\tau.
\end{align*}
We use the normalized versions of all B-splines with $n\leq7$, i.e., $\dicM$ contains all $96$ shifts of $B_n$, $n\leq7$. We have chosen a time period with low precipitation and thus the data can be sparsely approximated using splines. Note that a B-spline of order $n$ has a support of length $n$. Thus, the order directly correlates to the duration of the precipitation. Since a single precipitation (e.g., rain) will be registered at several stations, we have an underlying structure in the data.

Choose $\ptM_i$ to be the position coordinates of the $i$-th station and $\ptP_j\in\N^2$ as the shift and order of the corresponding spline. We use $d_{\domM}$ as the geodetic distance and set $\p=30$km. Let $d_{\domP}=\|\cdot\|_{\infty}$ and $\pp=1/15$, i.e., neither the duration nor the time of occurrence should change more than $2$ hours per $30$km. We perform $L=100$ iterations of GM-OMP.

\begin{figure}
\centering
\includegraphics[width=43mm]{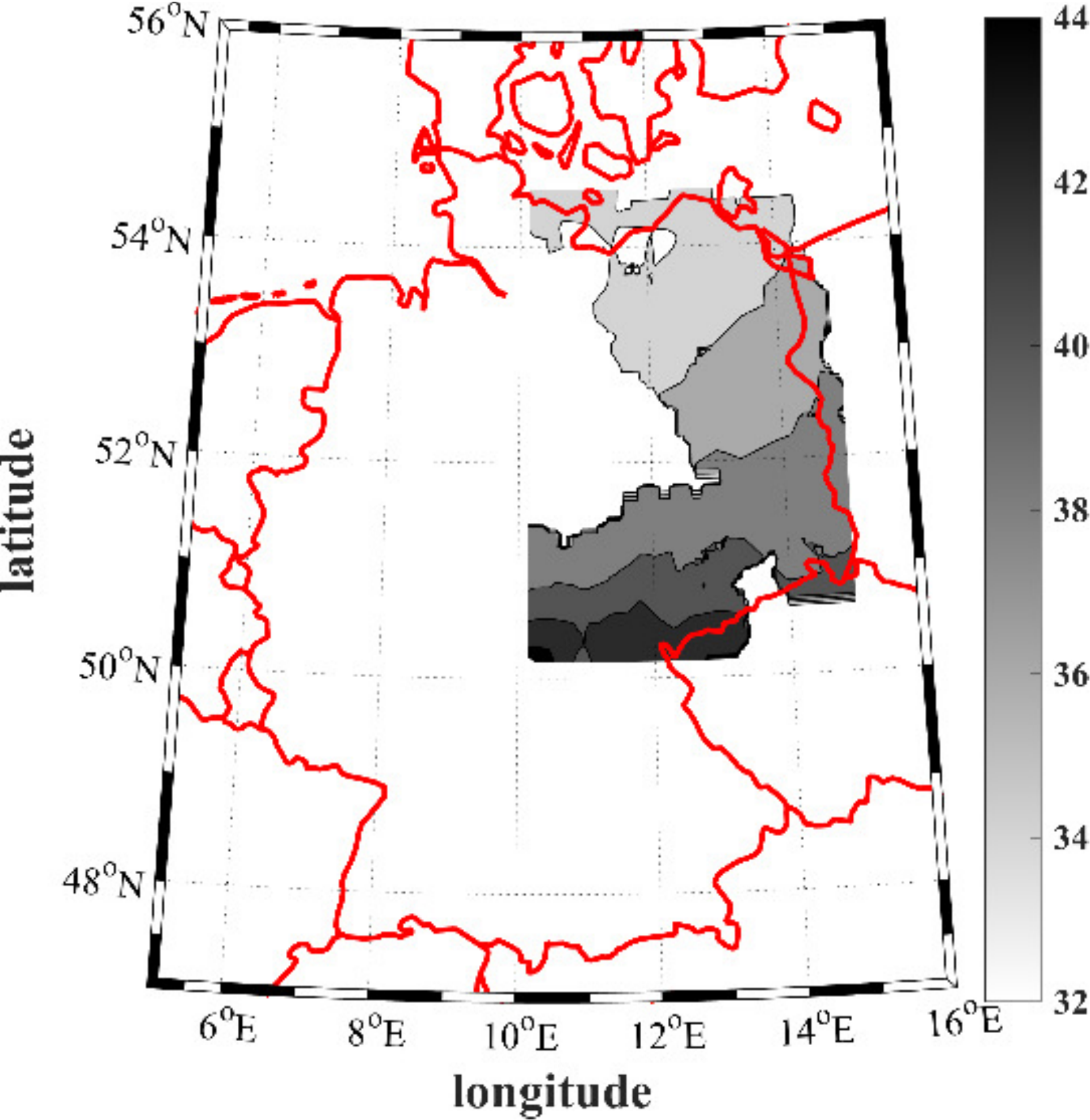}
\caption{Time of occurrence of the largest precipitation \cite{MMAP}.}
\label{fig:event70}
\end{figure}

\begin{figure}
\centering
\includegraphics[width=43mm]{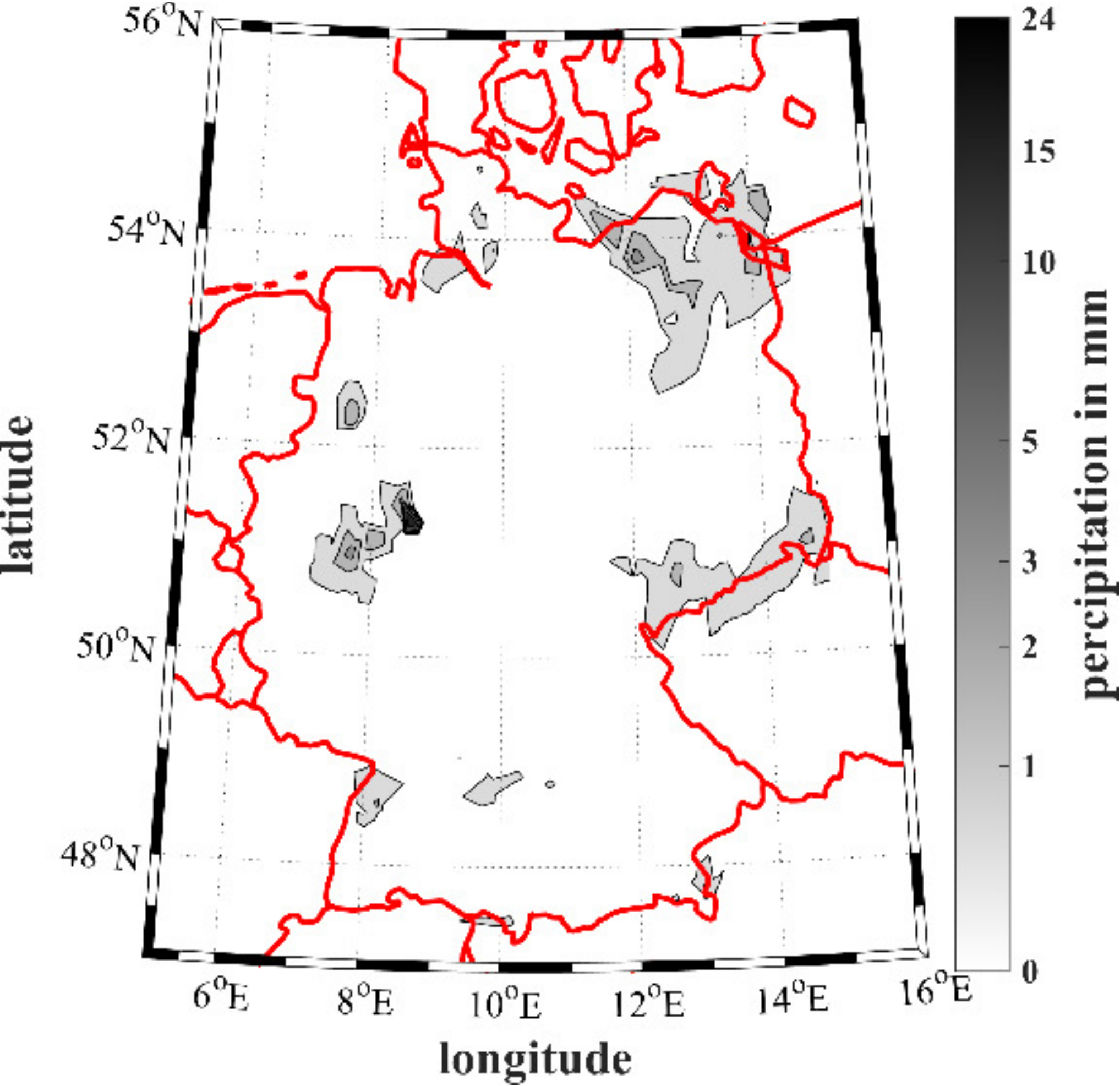}
\includegraphics[width=43mm]{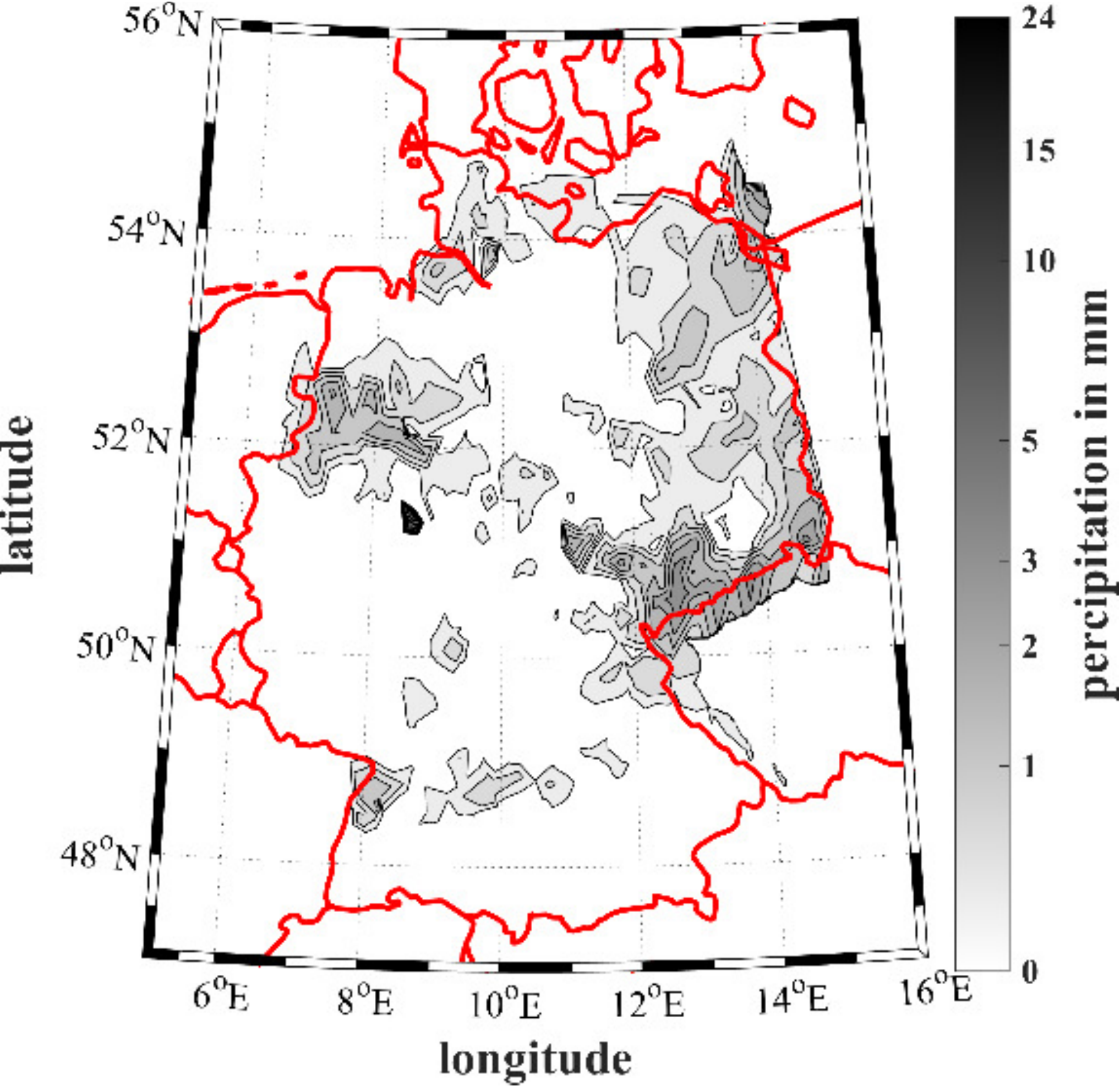}
\caption{Overall precipitation contour plot using only the $15$ strongest (left) or largest (right) events \cite{MMAP}.}
\label{fig:meanRek}
\end{figure}

Fig. \ref{fig:event70} shows the time of occurrence of the largest precipitation event. i.e., the structure that includes the most stations (here $156$). The mean duration is $1.28$h (mean B-spline order) and one can clearly recognize the event moving from north to south caused e.g., by wind. In Fig. \ref{fig:meanRek} we reconstructed the overall precipitation (see Fig. \ref{fig:metData}) using only $15$ structures. In the left figure we choose the first $15$ structure, i.e., those with the strongest precipitation; for the right figure the $15$ largest events were used. While the strongest events contain the precipitation peaks, the largest events can better reconstruct the overall structure from Fig. \ref{fig:meanRek}.

\section{Conclusion}

We presented a generalized orthogonal matching pursuit for multiple measurements. The algorithm is able to recognize and reconstruct more general sparsity pattern in the solution as other algorithms for multiple measurements. Moreover, GM-OMP allows efficient post processing for each pattern. These patterns can provide crucial information in application which was exemplary demonstrated in two practical examples. Two parameters allow an adaption of the feasible patterns to the application and make the algorithm more flexible. The advantages of GM-OMP were shown in comparison to other techniques and confirmed by first theoretical results.

\section*{Acknowledgements}

The author thanks Mannesmann Salzgitter GmbH for providing the ultrasonic data used in this paper. This work is supported by BMBF joined research project ZeMat (grant number: 05M13MGA).

\end{document}